\documentclass[centertags,12pt]{amsart}
\usepackage{latexsym}
\usepackage{amssymb}

\numberwithin{equation}{section}
\makeatletter
\renewcommand{\subsection}{\@startsection
{subsection}{2}{0mm}{\baselineskip}{-0.25cm}
{\normalfont\normalsize\bf}}
\makeatother

\newtheorem{theorem}{Theorem}[section]
\newtheorem{proposition}[theorem]{Proposition}
\newtheorem{lemma}[theorem]{Lemma}
\newtheorem{corollary}[theorem]{Corollary}

\newtheorem{result}[theorem]{Result}
{\theoremstyle{definition}
\newtheorem*{definition*}{Definition}

\newtheorem{remark}[theorem]{Remark}
\newtheorem*{proposition*}{Proposition}
\newtheorem*{corollary*}{Corollary}
\newtheorem*{lemma*}{Lemma}
\theoremstyle{remark}

\def\D{\mathbf D}

\def\P{\mathbf P}

\def\cA{\mathcal A}
\def\cB{\mathcal B}
\def\cC{\mathcal C}
\def\cD{\mathcal D}

\def\cF{\mathcal F}
\def\cG{\mathcal G}

\def\cH{\mathcal H}

\def\cO{\mathcal O}
\def\cP{\mathcal P}

\def\cS{\mathcal S}
\def\cU{\mathcal U}

\def\cX{\mathcal X}
\def\cY{\mathcal Y}

\def\Aut{\mbox{\rm Aut}}
\def\ord{\mbox{\rm ord}}
\def\deg{\mbox{\rm deg}}

\def\div{\mbox{\rm div}}
\def\min{{\rm min}}
\def\Aut{\mbox{\rm Aut}}
\def\Im{\mbox{\rm Im}}
\def\Ker{\mbox{\rm Ker}}

\def\Diff{\mbox{\rm Diff}}
\def\Alt{\mbox{\rm Alt}}
\def\Sym{\mbox{\rm Sym}}

\def\fq{{\mathbb F}_q}

\def\bfq{{\bar{\mathbb F}}_q}
\def\frx{{\mathbf Fr}}

\sloppy

\begin{document}
\author[M.~Giulietti]{Massimo Giulietti}
\author[G.~Korchm\'aros]{G\'abor Korchm\'aros}
\author[F.~Torres]{Fernando Torres}
\thanks{This research was supported by Italian Ministry MIUR, project {\em Strutture Geometriche, Combinatoria e loro applicazioni}, PRIN 2001-2002, and by GNSAGA. The third author was supported by the Grant SB2000-0225 from the ``Secretaria de Estado de Educaci\'on y Universidades del Ministerio de Educaci\'on, Cultura y Deportes de Espa\~na''}

\title[Quotient curves of the Deligne-Lusztig curve of Suzuki type]{Quotient curves of the Deligne-Lusztig curve\\
of Suzuki type}

\address{Dipartimento di Matematica e Informatica, Universit\`a di
Perugia, 06123 Perugia, Italy}
\email{giuliet@dipmat.unipg.it}
\address{Dipartimento di Matematica, Universit\`a della Basilicata, Campus Universitario Contrada Macchia Romana, 85100,
Potenza, Italy}
\email{korchmaros@unibas.it}
\address{IMECC-UNICAMP, Cx.P. 6065, Campinas, 13083-970, SP-Brazil}
\email{ftorres@ime.unicamp.br}
\address{Current Address: Dpto. Algebra, Geometr\'ia y
Topolog\'ia. Fac. Ciencias - Universidad de Valladolid, c/ Prado de la
Magdalena s/n 47005, Valladolid, Castilla, Spain}
\email{ftorres@agt.uva.es}

    \begin{abstract}
Inspired by a recent paper of Garcia, Stichtenoth and Xing [2000, {\em Compositio Math.} {\bf 120}, 137--170], we investigate the quotient curves of the Deligne-Lusztig curve associated to the Suzuki group $\cS z(q)$.
    \end{abstract}

\maketitle

    \section{Introduction}

The Deligne-Lusztig curve of Suzuki type (shortly DLS-curve) is
the (projective geometrically irreducible, non-singular) algebraic
curve defined to be the non-singular model over the finite field
$\mathbb{F}_q$ of the (absolutely irreducible) plane curve $\cC$
of equation $X^{q_0}(X^q+X)=Y^q+Y$, where $q_0=2^s,s\geq 1$ and
$q=2q_0^2$. Several authors have studied  the DLS-curve also in
connection with coding theory, see \cite{deligne-lusztig},
\cite{ft2}, \cite{hansen}, \cite{hansen-pedersen},
\cite{hansen-sti}, \cite{henn}. Here we only mention that the
DLS-curve has genus $g=q_0(q-1)$ and that the number of its
$\mathbb{F}_q$-rational points is $q^2+1$. Actually, the two
latter properties characterize the DLS-curve, see \cite{ft2}.
The automorphism group of the DLS-curve is the Suzuki group $\mathcal{S}z(q)$. In this paper, we investigate the quotient curves of the DLS-curve arising from the subgroups of ${\cS}z(q)$. For tame covering, that
is for subgroups of odd order, we obtain an exhaustive list of such
curves as given in the following theorem.
    \begin{theorem}\label{t1} Let $\cX$ be a tame quotient curve of the DLS-curve$.$ Then one of the following holds$.$
\begin{itemize}
\item[ I)] $r$ is any divisor of $q-1,$ $\cX$ has genus
$g=\frac{q-1}{r}q_0$ and is a non-singular model over $\fq$ of the plane curve of equation
   $$
Y^{(q-1)/r}\Big(
1+\sum_{i=0}^{s-1}X^{2^i(2q_0+1)-(q_0+1)}(1+X)^{2^i}\Big)
=(X^{q_0}+1)(Y^{2(q-1)/r}+X^{q-1})\, ,
   $$
\item[ II)] $r$ is any divisor of $q+2q_0+1,$ $\cX$ has genus
$g=\frac{q+2q_0+1}{r}(q_0-1) +1$ and is a non-singular model over
${{\mathbb F}_{q^4}}$ of the plane curve of equation
   $$Y^{(q+2q_0+1)/r}
\Big( 1+\sum_{i=0}^{s-1}
X^{2^iq_0}(1+X)^{2^i(q_0+1)-q_0}+X^{q/2}\Big)=
X^{q+2q_0+1}+Y^{2(q+2q_0+1)/r}\, ,
  $$
\item[ III)] $r$ is any divisor of $q-2q_0+1,$ $\cX$ has genus
$g=\frac{q-2q_0+1}{r}(q_0+1) -1$ and is a non-singular model over
${{\mathbb F}_{q^4}}$ of the plane curve of equation
  $$
bY^{(q-2q_0+1)/r}\Big(1+\sum_{i=0}^{s-1}X^{2^i(2q_0+1)-(q_0+1)}(1+X)^{2^i}\Big)=
(X^{q-2q_0+1}+Y^{2(q-2q_0+1)/r})(X^{q_0-1}+X^{2q_0-1})\, ,
  $$
\end{itemize}
where $b=\lambda^{q_0}+\lambda^{q_0-1}+\lambda^{-q_0}+\lambda^{-q_0+1}$
and $\lambda\in {{\mathbb F}_{q^4}}$ is an element of order $q-2q_0+1.$
    \end{theorem}
A similar complete list for non-tame coverings cannot be produced because the Suzuki group contains a huge number of pairwise non-isomorphic subgroups of even order. Our contribution consists in proving the existence of non-tame quotient curves of the DLS-curve of genus $g$ as given in Theorem \ref{t2}. For some of these curves we also provide a plane
equation, see Theorem \ref{t3}.
    \begin{theorem}\label{t2} Let $v,u,r$ be  positive integers$.$ For the following values of $g$ the DLS-curve has a quotient curve $\cX$ of genus $g.$
    \begin{itemize}
\item[ i)] $g=2^{s-u+v}(2^{2s+1-v}-1),$ $v\le 2s+1,$  $u\le v+\log_2{(v+1)},$
\item[ ii)] $g=\frac{1}{r}2^s(2^{2s+1-v}-1),$ $v\le 2s+1,$ $r|(q-1),$  $r|(2^{2s+1-v}-1),$
\item[ iii)] $g=\frac{q_0(q-r-1)}{2r},$ $r|(q-1),$
\item[ iv)] $ g=\frac{q_0(q-1)-1}{r}-(q_0-1),$ $r\mid (q+2q_0+1),$
\item[ v)] $ g=\frac{q_0(q-1)-1}{r}-(q_0-1),$ $r\mid (q-2q_0+1),$
\item[ vi)] $g=\frac{1}{2}\left[ \frac{q_0(q-1)-1}{r}- (q_0-1)\right],$ $r|(q+2q_0+1),$
\item[ vii)] $g=\frac{1}{2}\left[ \frac{q_0(q-1)+1}{r}- (q_0+1)\right],$ $r|(q-2q_0+1),$
\item[ viiii)] $g=\frac{1}{4}\left[ \frac{q_0(q-1)-1}{r}- (q_0-1)\right],$ $r|(q+2q_0+1),$
\item[ ix)] $g=\frac{1}{4}\left[ \frac{q_0(q-1)+1}{r}- (q_0+1)\right],$ $r|(q-2q_0+1),$
\item[ x)] $g=\frac{q_0(q-1)-1+({\bar{q}}^2+1){\bar{q}}^2({\bar{q}}-1)+\Delta}{({\bar{q}}^2+1){\bar{q}}^2({\bar{q}}-1)},$
${\bar q}=2^{2{\bar s}+1}$, ${\bar s}|s$ , $(2{\bar s}+1)|(2s+1),$
$\Delta:=(\bar{q}^2+1)[(2q_0+2)(\bar{q}-1)+2\bar{q}(\bar{q}-1)]+
\bar{q}^2(\bar{q}^2+1)(\bar{q}-2)+\bar{q}^2(\bar{q}+2\bar{q}_0+1)
(\bar{q}-1)(\bar{q}-2\bar{q}_0),$
\item[ xi)] $g=2^{4}(2^{9-v}-1),$ $3\le v\le 2s+1,$ for $q=512.$
    \end{itemize}
    \end{theorem}
    \begin{theorem}\label{t3} \begin{itemize}
\item[ i')] For $u=v,$ $v|(2s+1)$ a non-singular model over $\fq$ of the plane curve of equation
    $$
X^{2q_0}(X^q+X)=b\sum_{i=0}^{(2s+1/v)-1}Y^{(2^v)^i}
   $$
is a quotient curve of the DLS-curve of genus $g$ as in i)$.$
\item[ ii')] For $u=2,$ $v=1$ a non-singular model over $\fq$ of the plane curve of equation
   $$
\sum_{i=0}^{2s}X^{2^i}+\sum_{i=0}^{s}X^{2^i}\Big( \sum_{j=i}^s
X^{2^j}\Big) +\sum_{i=s+1}^{2s}
X^{2^i}\Big(\sum_{j=0}^{i-s-2}X^{2^j}\Big)^{2q_0}= \sum
_{i=0}^{2s}Y^{2^i}
  $$
is a quotient curve of the DLS-curve of genus $g$ as in i)$.$
\item[ iii')] A non-singular model over $\fq$ of the plane curve of equation
  $$
1+\sum_{i=0}^{s-1}X^{2^i(2q_0+1)-(q_0+1)}(1+X)^{2^i}=\sum
(-1)^{i+j}\frac{(i+j-1)!k}{i!j!}Y^i(X^{rj}(X^{q_0}+1)
  $$
where the summation is extended over all pairs $(i,j)$ of non--negative
integers with $i+2j=(q+2q_0+1)/r,$ is a quotient curve of the DLS-curve of genus $g$ as in iii)$.$
\item[ iv')] A non-singular model over ${\mathbb F}_{q^4}$ of the plane curve of equation
  $$Y^{(q+2q_0+1)/r}\Big(1+\sum_{i=0}^{s-1}
X^{2^iq_0}(1+X)^{2^i(q_0+1)-q_0}+X^{q/2}\Big)=
X^{q+2q_0+1}+Y^{2(q+2q_0+1)/r}
  $$
is a quotient curve of the DLS-curve of genus $g$ as in iv)$.$
\item[ v')] Let $b$ as in III)$.$ A non-singular model over ${\mathbb F}_{q^4}$ of the plane curve of equation
  $$
bY^{\frac{q-2q_0+1}{r}}\Big(1+\sum_{i=0}^{s-1}X^{2^i(2q_0+1)-q_0-1)}(1+X)^{2^i}\Big)
= (X^{q-2q_0+1}+Y^{\frac{2(q-2q_0+1)}{r}})(X^{q_0-1}+X^{2q_0-1})
  $$
is a quotient curve of the DLS-curve of genus $g$ as in v)$.$
\item[ vi')] A non-singular model over ${\mathbb F}_{q^4}$ of the plane curve of equation
   $$
1+\sum_{i=0}^{s-1} X^{2^iq_0}(1+X)^{2^i(q_0+1)-q_0}+X^{q/2}=\sum
(-1)^{i+j}\frac{(i+j-1)!}{i!j!}X^{ri}Y^j
   $$
is a quotient curve of the DLS-curve of genus $g$ as in vi)$.$
\item[ vii')] Let $b$ be as in III)$.$ A non-singular model over ${\mathbb F}_{q^4}$ of the plane curve of equation
   $$
b\Big(1+\sum_{i=0}^{s-1}X^{2^i(2q_0+1)-(q_0+1)}(1+X)^{2^i}\Big)
=(X^{q_0-1}+X^{2q_0-1})\sum
(-1)^{i+j}\frac{(i+j-1)!}{i!j!}X^{ri}Y^j
   $$
where the summation is extended over all pairs $(i,j)$ of non--negative integers with $i+2j=(q+2q_0+1)/r,$ is a quotient
curve of the DLS-curve of genus $g$ as in vii)$.$
   \end{itemize}
   \end{theorem}
A motivation for the present work comes from the current interest
in curves over finite fields with many rational points, see van der Geer's survey \cite{survey}. Indeed, the number of $\fq$-rational
points of a curve of genus $g$ which is $\fq$-covered by the
DLS-curve is $N=1+q+2q_0g$ (see Proposition \ref{nptscovered}) and
this value is in the interval from which the entries of the tables
of curves with many rational points are taken for $g\le 50$, $q\le
128$ in \cite{geer-vlugt}. Especially, both III) and v) for
$q=32$, $r=5$ provide new entries, namely $g=10$, $N=113$, and
$g=24$, $N=225$, and hence they imply $N_{32}(10)\ge 113$,
$N_{32}(24)\ge 225$. Furthermore, the following table shows those
values of $(q,g)$ for which some of the curves in the above
theorems attain the largest value of $\fq$-rational points for
which an $\fq$-rational curve of genus $g$ is previously known to
exist.
    \begin{center}
    \begin{tabular}{|c|c|l|c|}
\hline
$(q,g)$ & $N$ & Conditions & References\\
\hline\hline
$(8,6)$ & $33$ & i) with $v=1$, $u=1$ & \cite{sti2} \\
\hline
$(32,12)$ & $129$  & i) with $v=3$, $u=4$& \cite{vdg1} \\
\hline
$(32,28)$ & $257$ & i) with $v=2$, $u=3$ & \cite{vdg1} \\
\hline
$(32,30)$ & $273$ & i) with $v=1$, $u=2$ & \cite{vdg1} \\
\hline
$(128,8)$ & $257$ & I) with $r=127$, or & \cite{wirtz}\\
 &  & II) with $r=145$, or &  \\
 &  & III) with $r=113$, or & \\
 & & i) with $v=6$, $u=6$ & \\
\hline
$(128,14)$ & $353$ & i) with $v=4$, $u=6$, or & \cite{vdg2} \\
& & iv) with $r=29$ & \\
\hline
$(128,24)$ & $513$ & i) with $v=5$, $u=5$ & \cite{vdg1} \\
\hline
$(128,28)$ & $577$ & i) with $v=4$, $u=5$ & \cite{vdg1} \\
\hline
$(128,30)$ & $609$ & i) with $v=3$, $u=5$ & \cite{vdg2} \\
\hline
    \end{tabular}
    \end{center}

    \section{Preliminary results on the Suzuki group $\mathcal{S}z(q)$}\label{prel}

The Suzuki group ${\mathcal{S}}z(q)$ has been the
subject of numerous papers in finite geometry and permutation
group theory. Here, we only summarise those results on the
structure of ${\mathcal{S}}z(q)$ which play a role in the present
work. For more details, the reader is referred to \cite[Chapter
XI.3]{huppblack3}, \cite{Lun}, and \cite{tits}.
    \begin{result}\label{Ris1} Let $q_0=2^s,$ $s \geq 1,$ and $q=2q_0^2.$  For $a,c\in\fq,$ let
    $$
\widetilde{T}_{a,c}:= \left( \begin{array}{ccccc} 1 & 0 & 0 & 0
\\ a & 1 & 0 &  0 \\ c  & a^{2q_0}  & 1 & 0 \\
ac+a^{2q_0+2}+c^{2q_0} & a^{2q_0+1}+c & a & 1
\end{array} \right)\, .
   $$
The set $\widetilde{\bf T}=\{\widetilde{T}_{a,c}|a,c\in\fq\}$
is a group of exponent $4,$ class $2$ and order $q^2.$ For any
$b\in\fq$ with $b\neq 0,$ $\widetilde{\bf T}$ is isomorphic to the
$2$-dimensional linear group over $\fq$ consisting of all linear
transformation $(X,Y)\mapsto (X+a,b^{-1}a^{2q_0}X+Y+b^{-1}c)$ with
$a,c$ ranging over $\fq.$ In particular$,$ the center
$Z(\widetilde{\bf T})$ is an elementary abelian group of order
$2^{2s+1}$ whose non--trivial elements are those of order $2$ of
$\widetilde{\bf T}.$
    \end{result}
    \begin{result}\label{Ris2} For $d\in\fq$ and $d\neq 0,$ let
    $$
\widetilde{N}_d:=
 \left( \begin{array}{ccccc}
d^{-q_0-1} & 0 & 0 & 0\\ 0 & d^{-q_0} & 0 & 0\\ 0 & 0 & d^{q_0}&
0\\ 0 & 0 & 0 & d^{q_0+1}
\end{array} \right)\, ,
   $$
The set $\widetilde{\bf N}=\{\widetilde{N}_d|d\in\fq, d\neq
0\}$ is a group of order $q-1,$ isomorphic to the multiplicative
group of $\fq$. Furthermore$,$ $\widetilde{\bf N}$ normalises
$\widetilde{\bf T}$ and induces a fixed point free automorphism on
$\widetilde{\bf T}$ such that $\widetilde{\bf T}\widetilde{\bf N}$
is a Frobenius group with kernel $\widetilde{\bf T}.$
    \end{result}
Let
   $$
\widetilde{W}:= \left(
\begin{array}{ccccc} 0 & 0 & 0 & 1 \\ 0 & 0 & 1& 0\\ 0 & 1& 0 &
0\\ 1& 0 & 0 & 0
\end{array} \right)\, .
   $$
The Suzuki group ${\mathcal{S}}z(q)$ is defined as the
$4$-dimensional linear group generated by $\widetilde{\bf
T},\widetilde{\bf N}$ and $\widetilde{W}$, and the normaliser
$N_{{\mathcal{S}}z(q)}(\widetilde{\bf N})$ is the dihedral group
generated by $\widetilde{\bf N}$ together with $\widetilde{W}$.
Since ${\mathcal{S}}z(q)$ is a simple group, ${\mathcal{S}}z(q)$
faithfully induces a linear collineation group of the
$3$-dimensional projective space $\P^3(\fq)$. It turns out that
${\mathcal{S}}z(q)$ preserves the ovoid $\cO_3$ of $\P^3(\fq)$
consisting of the points $(0:0:0:1)$ and
$(1:x:y:xy+x^{2q_0+1}+y^{2q_0})$ with $x,y$ ranging over $\fq$.
More precisely, ${\mathcal{S}}z(q)$ is faithfully represented on
$\cO_3$ as a $2$-transitive permutation group of Zassenhaus type
of order $(q^2+1)q^2(q-1)$.
The stabiliser of the point $(0:0:0:1)$ under ${\mathcal{S}}z(q)$
is the Frobenius group $\widetilde{\bf T}\widetilde{\bf N}$.
Furthermore, ${\mathcal{S}}z(q)$ contains two conjugacy classes of
subgroups of Singer type, one consisting of cyclic subgroups
$\widetilde{\bf D}^+$ of order $q+2q_0+1$ and the other of cyclic
subgroups $\widetilde{\bf D}^-$ of order $q-2q_0+1$. The
normaliser $N_{{\mathcal{S}}z(q)}(\widetilde{\bf D}^+)$ has order
$4(q+2q_0+1)$ and is the semidirect product of $\widetilde{\bf
D}^+$ by  a cyclic group of order $4$. In particular,
$N_{{\mathcal{S}}z(q)}(\widetilde{\bf D}^+)$ is a Frobenius group
with Frobenius kernel ${\widetilde{\bf D}}^+$. All these results
hold true for ${\widetilde{\bf D}}^-$.
    \begin{result}\label{Ris3} The conjugates of the above subgroups$,$  namely $\widetilde{\bf T},\widetilde{\bf N},\widetilde{\bf D}^+,$ and
$\widetilde{\bf D}^-,$ form a partition of ${\mathcal{S}}z(q).$ If
$c(G)$ denotes the number of conjugates of a subgroup $G,$ then
$c(\widetilde{\bf T})=q^2+1,$ $c(\widetilde{\bf
N})=\frac{1}{2}q^2(q^2+1),$ $c(\widetilde{\bf
D}^+)=\frac{1}{4}q^2(q-2q_0+1)(q-1)$ and $c(\widetilde{\bf
D}^-)=\frac{1}{4}q^2(q+2q_0+1)(q-1).$
    \end{result}
In some cases, ${\mathcal{S}}z(q)$ contains subgroups isomorphic
to the Suzuki group over a subfield ${\mathbb F}_{\bar q}$ of $\fq$. This
occurs if and only if $\bar{q}=2^{2\bar{s}+1}$ with a divisor
$\bar{s}$ of $s$ such that $2\bar{s}+1$ divides $2s+1$. For such a
$\bar{q}$, there is just one conjugacy class in
${\mathcal{S}}z(q)$.
    \begin{result}\label{Ris4} Any subgroup of ${\mathcal{S}}z(q)$ is conjugate to either to ${\mathcal{S}}z(\bar{q}),$ or to a subgroup of one of
the following groups$:$ $\widetilde{\bf T}\widetilde{\bf N},$
$N_{{\mathcal{S}}z(q)}(\widetilde{\bf N}),$
$N_{{\mathcal{S}}z(q)}(\widetilde{\bf D}^+),$
$N_{{\mathcal{S}}z(q)}(\widetilde{\bf D}^-).$
    \end{result}
In studying ${\mathcal{S}}z(q)$ as an automorphism group of the
DLS-curve, we will need a suitable representation of
${\mathcal{S}}z(q)$ as a linear collineation group of $\P^4(\fq)$.
     \begin{result}\label{Ris5} For $a,c\in \fq,$ let
   $$
T_{a,c}:= \left(
\begin{array}{ccccc} 1 & 0 & 0 & 0 & 0
\\ a                      & 1 & 0 & 0 & 0
\\ a^{q_0+1}+c^{q_0}      & a^{q_0}      & 1 & 0 & 0
\\ c                        & a^{2q_0}     & 0 & 1 & 0 \\
ac+a^{2q_0+2}+c^{2q_0} & a^{2q_0+1}+c & 0 & a & 1
\end{array} \right)\, .
   $$
For $d\in\fq$ with $d\neq 0,$ let
   $$
N_d:=
 \left( \begin{array}{ccccc}
d^{-q_0-1} & 0 & 0 & 0 & 0\\ 0 & d^{-q_0} & 0 & 0& 0\\ 0 & 0 & 1 &
0 & 0\\ 0 & 0& 0 & d^{q_0}& 0\\
 0& 0 & 0 & 0 & d^{q_0+1}
\end{array} \right)\, , \qquad  W:= \left(
\begin{array}{ccccc} 0 & 0 & 0 & 0 & 1
\\
 0 & 0 & 0 & 1& 0\\
0 & 0 & 1& 0 & 0\\ 0 & 1& 0 & 0 & 0\\ 1& 0 & 0 & 0 & 0
\end{array} \right)\, .
     $$
Let $[ T_{a,c}],$ $[ N_d],$ and $[W]$ be the linear collineations
of $\P^4(\fq)$ associated to $T_{a,c},$ $N_d,$ and $W,$
respectively. Then the group $[\cG]$ generated by them is
isomorphic to ${\mathcal{S}}z(q).$ Let
   $$
\cO_4=\{(1:u:u^{q_0+1}+v^{q_0}y^{q_0}:bv:buv+u^{2q_0+2}+b^{2q_0}v^{2q_0})|u,v\in\fq\}
\cup\{(0:0:0:0:1)\}\, .
   $$
Then $[\cG]$ preserves $\cO_4$ and acts on
it as ${\cS}z(q)$ in its unique $2$-transitive permutation
representation$.$ The full collineation group of $\P^4(\fq)$
preserving $\cO_4$ is isomorphic to $\Aut({\cS}z(q))$ and is the
semidirect product of $[\cG]$ by the non-linear group $[\cF]$
of order $2s+1$ of $\P^4(\fq)$ generated by the collineation
$(X_0:X_1:X_2:X_3:X_4) \mapsto (X_0^2:X_1^2:X^2_2:X_3^2:X_4^2)$.
   \end{result}
   \begin{proof} Let $[\cG]$ be the linear collineation group of $\P^4(\fq)$
generated by $[T_{a,c}],[N_d]$ and $[W]$, with $a,c,d$ ranging in
$f_q$ and $d\neq 0$. Then $[\cG]$ preserves the set $\cO_4$
consisting of the points
$(1:x:x^{q_0+1}+b^{q_0}y^{q_0}:by:bxy+x^{2q_0+2}+b^{2q_0}y^{2q_0})$
and $(0:0:0:0:1)$. Following the method used for the explicit
construction of ${\cS} z(q)$ in \cite[Chapter IV]{Lun} one can
extend Tits' result to $[\cG]$, see \cite[Theorem 21.8]{Lun}:
$[\cG]$ acts doubly transitively on $\cO_4$; if $[G]\in[\cG]$,
then either $[G]=[N_d][T_{a,c}]$ for exactly one triple $(d,a,c)$
with $d(\neq 0),a,c\in \fq$, or $[G]=[N_d][T_{a,c}][W][T_{e,f}]$,
for exactly one quintuple $(d,a,c,e,f)\in \fq$, with $d(\neq
0),a,c,e,f\in \fq$. Also, Theorem 21.11 and the Suzuki Tits
Theorem 22.6 in \cite{Lun} hold true: $[\cG]$ is a simple group of
order $q^2(q^2+1)(q-1)$. This implies that $[\cG]\cong {\cS}z(q)$.
We remark that the image of $[G]\in[\cG]$ under such an
isomorphism can be easily obtained by deleting the third row and
the third column of $G$. Since the non-linear collineation $F:
(X_0:X_1:X_2:X_3:X_4) \mapsto (X_0^2:X_1^2:X^2_2:X_3^2:X_4^2)$
preserves $\cO_4$, we have to prove that every collineation $h$ of
$\P^4(\fq)$ preserving $\cO_4$ is in $[\cH][\cF]\cong
\Aut({\cS}z(q))$. If $h$ fixes $\cO_4$ pointwise, then $h$ is the
identity because $\cO_4$ is not contained in any proper subspace
of $\P^4(\fq)$. Thus $h$ acts on $\cO_4$ faithfully. On the other
hand, from the classification of finite $2$-transitive permutation
groups it follows that there are only three permutation groups of
degree $q^2+1$ containing ${\cS}z(q)$, namely $\Aut({\cS}z(q))$,
$\Alt_{q^2+1}$ and $\Sym_{q^2+1}$. It remains to show that the two
latter cases cannot actually occur in our situation. Take six
points on $\cO_4$ no five of them linearly dependent in
$\P^4(\fq)$. By the fundamental theorem of projective geometry,
the group generated by $F$ is the collineation group of
$\P^4(\fq)$ fixing each of these six points. On the other hand,
the stabilizer of six objects in $\Alt_n$ with $n \geq 10$ is not
a cyclic group. This proves the assertion.
    \end{proof}
    \begin{remark} From now on we will use the term of Suzuki group and the symbol ${\cS}z(q)$ to denote the four dimensional projective linear group
introduced in Result \ref{Ris5}. In particular, $\bf T,\bf N, \bf
D^+, \bf D^-$ will stand for the corresponding subgroups
$\widetilde{\bf T},\widetilde{\bf N}, \widetilde{\bf
D}^+,\widetilde{\bf D}^-$ under the isomorphism stated in Result
\ref{Ris5}.
    \end{remark}

    \section{Preliminary results on the DLS-curve}\label{dls}

Throughout the present paper $\cX$ will stand for the
DLS-curve over $\fq$. As we have mentioned in the Introduction
$\cX$ has genus $q_0(q-1)$ and contains exactly $q^2+1$
$\fq$-rational points. By the Serre-Weil explicit formulae (see
\cite{serre}, \cite{hansen}), the characteristic polynomial
$h_{\cX}(t)$ of the Frobenius morphism over $\fq$ on the Jacobian
variety of $\cX$ is $(t^2+2q_0t+q)^{q_0(q-1)}$. Given a curve $\cY
$ of genus $g$ which is $\fq$-covered by $\cX$, Lachaud's theorem
\cite{lachaud} implies that $h_{\cY}(t)=(t^2+2q_0t+q)^g$. By
\cite[V.1]{sti}, the following proposition follows.
   \begin{proposition}\label{nptscovered}
Let $\cY$ be a curve which is $\fq$-covered by $\cX$. Then the
number of $\fq$-rational points of $\cY$ is equal to $1+q+2gq_0,$
where $g$ is the genus of $\cY$.
   \end{proposition}
   \begin{remark} Let $\cY$ be a curve as in the previous proposition. One can show then that $\#\cY(\mathbb F_{q^2})=1+q^2, \#\cY(\mathbb F_{q^3})=1+q^3-2gq_0q$ and $\#\cY(\mathbb F_{q^4})=1+q^4+2gq^2$; i.e., $\cY$ is maximal over $\mathbb F_{q^4}$. As a matter of fact, the examples obtained so far in this paper give new insights toward the computation of the spectrum of genera of maximal curves over finite fields of characteristic two (compare with the examples in \cite{g-sti-x}, \cite{ckt1}, \cite{ckt2} and \cite{aq}).
   \end{remark}
The proposition below will be useful in the sequel.
    \begin{proposition}\label{functionfield}
For any $b\in \fq$, $b\neq 0$, there are elements $x,y\in \fq(\cX)$ such that
   $$\fq(\cX)=\fq(x,y), \quad x^{2q_0}(x^q+x)=b(y^q+y).$$
   \end{proposition}
   \begin{proof} We have $\fq(\cX)=\fq(x,t)$ with $x^{q_0}(x^q+x)=t^q+t$. Let $y=b^{-1}(x^{2q_0+1}+t^{2q_0})$, that is
$t^q=b^{q_0}y^q_0+x^{q+q_0}.$ Then
$\fq(\cX)=\fq(x,y)$. Furthermore,
$y^{q_0}=b^{-q_0}(x^{q+q_0}+t^q)=b^{-q_0}(x^{q_0+1}+t)$, and hence
$y^q=b^{-1}(x^{q+q_0}+t^{2q_0})$. Now, since
$y^q+y=b^{-1}(x^{q+q_0}+t^{2q_0}+x^{2q_0+1}+t^{2q_0})=
b^{-1}x^{2q_0}(x^q+x),$
the claim follows.
   \end{proof}
Let $\cC_b$ be the plane curve of equation
$X^{2q_0}(X^q+X)=b(Y^q+Y)$. $\cC_b$  has only one singular point,
namely the infinite point $Y_{\infty}$ of the $Y$-axis which point
is a $q_0$-fold point. We know from \cite{hansen-sti} that
$\bfq(\cX)$ has just one place centered at $Y_{\infty}$. Let
$P_\infty$ denote the corresponding point of $\cX$. From now on,
we fix a projective frame $A_0A_1A_2A_3A_4U$ in
$\P^4({\bar{\mathbb{F}}}_q)$ with fundamental vertices
$A_0:=(1:0:0:0:0),\ldots, A_4:=(0:0:0:0:1)$, and $U:=(1:1:1:1:1)$.
With the notation of Proposition \ref{functionfield}, let $f$ be
the morphism $f:\cX \rightarrow \P^4({\bar{\mathbb{F}}}_q)$ with
coordinate functions
 $$
f:=(f_0:f_1:f_2:f_3:f_4)\,, $$ such that $f_0:=1,\,f_1:=x,\,
f_2:=x^{q_0+1}+b^{q_0}y^{q_0}, \, f_3:=by,\,
f_4:=bxy+x^{2q_0+2}+b^{2q_0}y^{2q_0}$. They are uniquely
determined by $f$ up to a proportionality factor in $\bfq(\cX)$.
For each point $P\in\cX$, we have
$f(P)=((t^{-e_P}f_0)(P),\ldots,t^{-e_P}f_4)(P))$ where
$e_P=-{\rm min}\{v_P(f_0),\ldots,v_P(f_4)\}$ for a local parameter $t$
of $\cX$ at $P$. It turns out that $f(\cX)$ is a parametrised
curve not contained in any hyperplane of $\P^4(\bfq)$. For a point
$P\in f(\cX)$, the intersection multiplicity of $f(\cX)$ with a
hyperplane $H$ of equation $a_0X_0+\ldots+a_4X_4=0$ is
$v_P(a_0f_0+\ldots+a_4f_4)+e_P$, and the intersection divisor
$f^{-1}(H)$ cut out on $f(\cX)$ by $H$ is defined to be
$f^{-1}(H)=\div(a_0f_0+\ldots+a_4f_4)+E$ with $E=\sum e_pP$. We
have $v_{P_\infty}(f_1)=-q$, $v_{P_\infty}(f_3)=-2q_0-q$,
$v_{P_\infty}(f_2)=-q_0-q$, $v_{P_\infty}(f_4)=-2q_0-q-1$, see
Section 6 . Then $e_{P_\infty}=q+2q_0+1$, and the representative
$(f_0/f_4:f_1/f_4:f_2/f_4:f_3/f_4:1)$ of $f$ is defined on
$P_{\infty}$. Hence $f(P_{\infty})=A_4$, and $f_3/f_4$ is a local
parameter at $P_\infty$.
For a point $P\in\cX$, an integer $j$ is called a hermitian
$P$-invariant (cf. \cite{sv}) if there exists a hyperplane intersecting $f(\cX)$ at $f(P)$ with multiplicity $j$. There are exactly five pairwise
distinct hermitian $P$-invariants. Such integers arranged in
increasing order define the order sequence of $\cX$ at $P$. By
\cite{ft2}, the order sequence of $\cX$ at a point $P\in\cX$ is
either $(0,1,q_0+1,2q_0+1,q+2q_0+1)$ or $(0,1,q_0,2q_0,q)$
according as $P\in\cX(\mathbb{F}_q)$ or
$P\not\in\cX(\mathbb{F}_q)$. The linear system
   $$\{f^{-1}(H) \
| \ H \ \text{hyperplane in} \ P^4(\bfq)\}
   $$
is $|(q+2q_0+1)P_0|$ for $P_0\in \cX(\fq)$. Also, $(q+2q_0+1)P\sim
qP+2q_0\frx(P)+\frx^2(P)$ for every $P\in\cX$, where $\frx$ is the Frobenius morphism over $\fq$; see \cite{ft2}.
    \begin{proposition}\label{closed} $f(\cX)$ is a non-singular model defined over $\fq$ of the DLS curve.
     \end{proposition}
\begin{proof} We show that $f$ is a closed embedding. By the above discussion,
$f$ is bijective and $f(\cX)$ has no singular point.
\end{proof}
According to Proposition \ref{closed}, we will identify $f(\cX)$
with $\cX$. 
\begin{proposition}
\label{suz} The automorphism group $\Aut(\cX)$ of $\cX$ is
isomorphic to ${\cS}z(q)$ and acts on $\cX(\fq)$ as ${\cS}z(q)$ in
its unique $2$-transitive permutation representation.
\end{proposition}
\begin{proof}
For $a,c,d \in \fq$ with $d\neq 0$, we define the following
automorphisms of $\fq(\cX)$:
\begin{equation}
\label{gamma}
 \psi_{a,c} \ := \left \{ \begin{array}{lll}
         x \mapsto x+a, \\
         y \mapsto b^{-1}a^{2q_0}x+y+b^{-1}c;
                         \end{array}
          \right.
\qquad
 \gamma_d \ := \left \{ \begin{array}{lll}
         x \mapsto dx, \\
         y \mapsto d^{2q_0+1}y;
                         \end{array}
          \right.
\end{equation}
for $h:=bxy+x^{2q_0+2}+b^{2q_0}y^{2q_0}$,
\begin{equation}
\label{varphi}
\varphi \ := \left \{ \begin{array}{lll}
         x \mapsto by/h, \\
         y \mapsto bx/h.
                         \end{array}
          \right.
\end{equation}
Let $\Gamma$ be the automorphism group of $\bfq(\cX)$ generated by
$\psi_{a,c}$, $\gamma_d$ and $\varphi$. By straightforward
computations, $[W]f=f\varphi$, $[N_d]f=f\gamma_d$ and
$[T_{a,c}]f=f\psi_{a,c}$. This shows that there is a homomorphism
$\Gamma \mapsto {\cS} z(q)$. Actually, this homomorphism is an
isomorphism because the identity is the only automorphism of
$\bfq(\cX)$ which acts as the identity map on the set of all
places of $\bfq(\cX)$, or, equivalently, on the set of all points
of $\cX$. Result \ref{Ris5} yields
$\Aut_{\fq}(\cX)\cong{\cS}z(q)$. Thus, $\Gamma=\Aut_{\fq}(\cX)$.
Finally, $\Aut_{\fq}(\cX)=\Aut(\cX)$ by \cite{henn}.
\end{proof}
\begin{remark}
\label{remark2} For the rest of the paper, $\cX=f(\cX)$ is
chosen for a non-singular model over $\fq$ of the DLS-curve. Then
$\Aut(\cX)\cong{\cS}z(q)$, and $\Aut(\cX)$  acts on the
set of places of $\bfq(\cX)$ as ${\cS}z(q)$ on the
set of points of $\cX$. In particular, $\psi_{a,c}$, $\gamma_d$,
and $\varphi$ correspond to $[T_{a,c}]$, $[N_d]$ and $[W]$ under
such an isomorphism.
\end{remark}

    \section{Quotient curves arising from subgroups of
a cyclic subgroup of ${\cS} z(q)$ of order $q-1$}\label{ordq-1}

For a divisor $r>1$ of $q-1$, let $\cU$ be a subgroup of
${\cS}z(q)$ of order $r$. Up to conjugacy in ${\cS}z(q)$, we have
$\cU=\{N_d|d^r=1\}$. It is straightforward to check that $\cU$ has
exactly two fixed points on $\cX$, namely $A_0$ and $A_4$. Let
$\cX_{\cU}$ denote the quotient curve of $\cX$ associated to
$\cU$, and let $g_{\cU}$ be its genus. Since $\cU$ is a tame
subgroup, the Hurwitz genus formula gives
$2q_0(q-1)-2=r(2g_{\cU}-2)+2(r-1)$, whence
$g_{\cU}=\frac{1}{r}q_0(q-1).$
To find an explicit equation of $\cX_{\cU}$ we first determine a
plane (singular) model of $\cX$ on which not only $\bf N$ but also
$[W]$ acts linearly.
For every non-zero element $d$ in $\fq$, both $[N_d]$ and $[W]$
preserve the line $\ell$ joining $A_0$ and $A_4$, as well as the
plane $\alpha$ spanned by the other three fundamental points,
namely $A_1,A_2$, and $A_3$. Now, project $\cX$ from $\ell$ to
$\alpha$. The associated morphism $\cX \to \P^2(\bfq)$ is
$\pi_{\ell}: (1:x:x^{q_0+1}+b^{q_0}y^{q_0}:by:bxy+x^{2q_0+2}+b^{2q_0}y^{2q_0})
\mapsto(x:x^{q_0+1}+b^{q_0}y^{q_0}:by)$. In terms of linear
systems, $\pi_{\ell}$ is associated to the $2$-dimensional linear
series cut out on $\cX$ by hyperplanes through $\ell$. We will
need some computational results.
\begin{equation}
\label{vPx}
 v_{P}(x) \ = \left \{ \begin{array}{lll}
 1 & \mbox{\rm for\ $P=(1:0:(bc)^{q_0}:bc:(bc)^{2q_0})$ with
$c\in\fq\setminus \{0\}$\/}, \\
   1  & \mbox{\rm for\ $P=A_0$\/}, \\
  -q  & \mbox{\rm for\ $P=A_4$\/}, \\
  0    & \mbox{\rm otherwise.\/}
                         \end{array}
          \right.
\end{equation}
\begin{equation}
\label{vPy}
 v_{P}(by) \ = \left \{ \begin{array}{lll}
 1 & \mbox{\rm for\ $P=(1:a:a^{q_0+1}:0:a^{2q_0+2})$ with $a\in\fq\setminus \{0\}$\/}, \\
 2q_0+1  & \mbox{\rm for\ $P=A_0$\/}, \\
 -2q_0-q  & \mbox{\rm for\ $P=A_4$\/}, \\
  0    & \mbox{\rm otherwise.\/}
                         \end{array}
          \right.
\end{equation}
\begin{equation}
\label{vxq0+1+yq0} v_{P}(x^{q_0+1}+(by)^{q_0}) = \left \{
\begin{array}{lll}
 1 & \mbox{\rm for\ $P=(1:a:0:a^{2q_0+1}:a^{2q_0+2}),a\in\fq\setminus \{0\}$\/}, \\
 q_0+1  & \mbox{\rm for\ $P=A_0$\/}, \\
 -q_0-q    & \mbox{\rm for \ $P=A_4$\/}, \\
 0    & \mbox{\rm otherwise.\/}
                         \end{array}
          \right.
\end{equation}
\begin{equation}
\label{vPx+by}
 v_{P}(x+by) \ = \left \{ \begin{array}{lll}
 1 & \mbox{\rm for $P=(1:a:a^{q_0+1}+a^{q_0}:a:a^{2q_0+2})$ with $a\in\fq\setminus \{0,1\}$\/}, \\
2q_0+1 & \mbox {\rm for $P=(1:1:0:1:1)$\/}, \\
 1 & \mbox{\rm for $P=A_0$\/}, \\
-2q_0-q  & \mbox{\rm for\ $P=A_4$\/}, \\
 0    & \mbox{\rm otherwise.\/}
                         \end{array}
          \right.
\end{equation}
For $e_{P}^{\prime}:=-\min
\{v_P(x),v_P(x^{q_0+1}+(by)^{q_0}),v_P(by)\}$,
\begin{equation}
\label{eprimeP}
 e_{P}^{\prime} \ = \left \{ \begin{array}{lll}
 -1 & \mbox{\rm for\ $A_0$\/}, \\
 q+2q_0  & \mbox{\rm for\ $P=A_4$\/}, \\
  0    & \mbox{\rm otherwise.\/}
                         \end{array}
          \right.
\end{equation}
\begin{proof}  Notice that $t:=x$ is a local parameter
at $A_0$. From $x^{2q_0}(x^q+x)=by^q+by$, we have
$y=b^{-1}t^{2q_0+1}+b^{-1}t^{q+2q_0}+\ldots$, and hence
$x+by=t+t^{2q_0+1}+\ldots$, and
$x^{q_0+1}+(by)^{q_0}=t^{q_0+1}+t^{q+q_0}+\ldots$ whence the above
results for $A_0$ follow. Also, $v_{A_0}(h)=q+2q_0+1$ where
$h=bxy+x^{2q_0+2}+(by)^{2q_0}$. Since the involutory automorphism
$[W]$ of $\cX$ changes $A_0$ with $A_4$, and acts on
$\bar{\mathbf{F}}_q(\cX)$ as (\ref{varphi}), we have
$v_{A_4}(x)=v_{A_0}(by/h)=-q$,
$v_{A_4}(by)=v_{A_0}(bx/h)=-2q_0-q$,
$v_{A_4}(x^{q_0+1}+(by)^{q_0})=v_{A_0}((by/h)^{q_0+1}+(bx/h)^{q_0})=-q_0-q$,
and
$v_{A_4}(x+by)=v_{A_0}(by/h+x/h)=v_{A_0}(x+by)-v_{A_0}(h)=-2q_0-q$.So, the above results for $A_4$ hold.
 Now, we assume $A_0\neq P \neq A_4$, that is
$P=(1:a:a^{q_0+1}+(bc)^{q_0}:bc:abc+a^{2q_0+2}+(bc)^{2q_0})$ with
$a^{2q_0+q}+a^{2q_0+1}=(bc)^q+bc$ and either $a\neq 0$, or $a=0$
and $c\in\fq\setminus\{0\}$. In the latter case, $t:=x$ is a local
parameter and the local expansions of the coordinate functions of
$\pi_{\ell}$ are
$$x=t,\quad x^{q_0+1}+b^{q_0}y^{q_0}=(bc)^{q_0}+t^{q_0+1}+\ldots,
\quad by=bc+t^{2q_0+1}+\ldots.$$ Hence $v_P(x)=1$, but
$v_P(by)=0$. Also, $v_P(x^{q_0+1}+(by)^{q_0})=0=v_P(x+by)$. We
have to investigate the case $a\neq 0$. Since $t:=x+a$ is a local
parameter, we have $x=a+t, y=c+\alpha t+\beta t^i+\ldots$ for some
$i>1$ and $\beta\neq 0$. From $x^{2q_0}(x^q+x)=(by)^{q}+by$ we
deduce that $\alpha\neq 0$ and that either $i=2q_0+1$ or $i=2q_0$
according as $a$ belongs to $\fq$ or does not. Hence the local
expansions of the coordinate functions of $\pi_{\ell}$ are
$$x=a+t,\quad x^{q_0+1}+b^{q_0}y^{q_0}=a^{q_0+1}+(bc)^{q_0}+a^{q_0}t+\ldots,
\quad by=bc+a^{2q_0}t+t^{2q_0+1}+ \ldots$$ for $a
\in\fq\setminus\{0\}$, and
$$x=a+t,\quad x^{q_0+1}+b^{q_0}y^{q_0}=a^{q_0+1}+(bc)^{q_0}+a^{q_0}t+\ldots,
\quad by=bc+a^{2q_0}t+(a^q+a)t^{2q_0}+\ldots$$ for $a\not\in\fq.$
Thus $v_P(x)=0$. Furthermore, $v_P(by)\neq 0$ if and only if
$c=0$. More precisely, this only occurs when
$a\in\fq\setminus\{0\}$, $P=(a:a^{q_0+1}:a^{q_0}:0:a^{2q_0+2})$,
and $v_P(by)=1$. Also, $v_P(x^{q_0+1}+y^{q_0})\neq 0$ if and only
if $a^{q_0+1}+(cb)^{q_0}=0$. This condition is only satisfied by
$a\in\fq$. In fact, $a^{q_0+1}=(cb)^{q_0}$ together with
$a^{2q_0}(a^q+a)=(cb)^q+cb$ implies $cb=a^{2q_0+1}$ and hence
$a^q=a$. For $a\in\fq\setminus\{0\}$, we have
$P=(1:a:0:a^{2q_0+1}:a^{2q_0+2})$ and $v_P(x^{q_0+1}+y^{q_0})=1$.
Moreover, $v_P(x+by)\neq 0$ if and only if $a+bc=0$, that is $a\in
\fq\setminus\{0\}$ and $P:=(1:a:a^{q_0+1}+a^{q_0}:a:a^{2q_0+2})$.
More precisely, either $v_P(x+by)=1$, or $v_P(x+by)=2q_0+1$
according as $a\in\fq\setminus\{0,1\}$ or $a=1$. Finally, the
formula for $e_P^{\prime}$ follows from (\ref{vPx}) and
(\ref{vPy}) together with (\ref{vxq0+1+yq0}).
 \end{proof}
The homogeneous coordinates
$(X_1^{\prime}:X_2^{\prime}:X_3^{\prime})$ provide a natural
projective frame in $\alpha$ with fundamental triangle
$A_1A_2A_3$.
    \begin{lemma}\label{cyclicq-1} The plane curve $\pi_{\ell}(\cX)$ is
birationally $\fq$-isomorphic to $\cX$ and it has degree
$q+2q_0-1$. The action of $[N_d]$ on $\pi_{\ell}(\cX)$ is induced
by the linear automorphism
$(X_1^{\prime}:X_2^{\prime}:X_3^{\prime})\mapsto(d^{-q_0}
X_1^{\prime}:X_2^{\prime}:d^{q_0}X_3^{\prime})$, whereas the
action of $[W]$ by $(X_1^{\prime}:X_2^{\prime}:X_3^{\prime})$
$\mapsto(X_3^{\prime}:X_2^{\prime}:X_1^{\prime})$.
     \end{lemma}
     \begin{proof} Let $P:=(1:1:1:0:1)$,
$Q:=(1:u:u^{q_0+1}+b^{q_0}v^{q_0}:bv:buv+u^{2q_0+2}+b^{2q_0}v^{2q_0})$
be two points on $\cX$ such that $\pi_{\ell}(Q)=\pi_{\ell}(P)$. A
straightforward computation yields $Q=P$. Moreover,
$by/(x^{q_0+1}+b^{q_0}y^{q_0})$ is a local parameter at $P$ as
$v_P(x)=0$ and $v_P(y)=1$. This implies that $\pi_{\ell}$ is
birational. The second assertion follows from (\ref{eprimeP}) by
virtue of $\deg (\pi_{\ell}(\cX))=\sum e^{\prime}_P$. The third
assertion is easily deduced from the matrix representations of
$[N_d]$ and $[W]$.
     \end{proof}
To write an equation of $\pi_{\ell}(\cX)$ we will use in $\alpha$
the affine frame $(X^{\prime},Y^{\prime})$ arising from the above
projective frame $(X_1^{\prime}:X_2^{\prime}:X_3^{\prime})$ by
$X^{\prime}=X_1^{\prime}/X_2^{\prime}$, and
$Y^{\prime}=X_3^{\prime}/X_2^{\prime}$. Then $\pi_{\ell}(\cX)$ has
equation $F(X^{\prime},Y^{\prime})=0$ where
$F(X^{\prime},Y^{\prime})$ is an absolutely irreducible polynomial
with coefficients in $\fq$ which satisfies  $F(\xi,\eta)=0$ where
$\xi,\eta$ are defined to be
$\xi:=x/(x^{q_0+1}+b^{q_0}y^{q_0}),\eta:=by/(x^{q_0+1}+b^{q_0}y^{q_0})$.
      \begin{proposition}\label{propo-1}
The equation of $\pi_{\ell}(\cX)$ can be written in the form $$
G_1(X^{\prime}Y^{\prime})=G_2(X^{\prime}Y^{\prime})
({X^{\prime}}^{q-1}+{Y^{\prime}}^{q-1})\,, $$ with
$G_1(T),G_2(T)\in {\mathbb F}_{q}[T]$, $\deg (G_1(T))=
\frac{1}{2}q-1$ and $\deg (G_2(T))= q_0$.
    \end{proposition}
    \begin{proof} Let $F(X^{\prime},Y^{\prime})=\sum
a_{ij}{X^{\prime}}^i{Y^{\prime}}^j.$ Let $r_{\infty}$ be the line
of $\alpha$ with equation $X_2^{\prime}=0$.  By (\ref{vxq0+1+yq0})
and (\ref{eprimeP}), the intersection divisor of $r_{\infty}$ is
given by $\pi_{\ell}^{-1}(r_{\infty})=q_0A_0+q_0A_4+\sum
P(1:a:0:a^{2q_0+1}:a^{2q_0+2})$ with $a$ ranging over
$\fq\setminus\{0\}$. In other words, for the intersection
multiplicity of $\pi_{\ell}(\cX)$ with the line $r_{\infty}$ we have
   $$
I(\pi_{\ell}(\cX),r_{\infty};\pi_\ell(P)) = \left \{
\begin{array}{lll}
 1 & \mbox{\em for\ $P=(1:a:0:a^{2q_0+1}:a^{2q_0+2})$
 with  $a\in\fq\setminus \{0\}$\/}, \\
 q_0  & \mbox{\em for\ $P=A_0$ and $P=A_4$ \/}, \\
  0    & \mbox{\em otherwise.\/}
                         \end{array}
          \right.
   $$
This implies that the  monomials of degree $q+2q_0-1$ of
$F(X',Y')$ are ${Y'}^{q_0}{X'}^{q+q_0-1}$ and
${X'}^{q_0}{Y'}^{q+q_0-1}$. By the last
two claims in Lemma \ref{cyclicq-1}, there exist $e\, , e'\, \in
{\bar \fq}$ such that
   $$ \sum
a_{ij}{X^{\prime}}^i{Y^{\prime}}^j=e\sum a_{ij}\omega^{q_0(i-j)}
{X^{\prime}}^i{Y^{\prime}}^j\, ,
   $$
   $$ \sum
a_{ij}{X^{\prime}}^i{Y^{\prime}}^j=e'\sum a_{ij}
{Y^{\prime}}^i{X^{\prime}}^j\, ,
   $$
where $\omega$ is a generator
of the multiplicative group of $\fq$.  Then
$a_{ij}=ea_{ij}\omega^{q_0(i-j)},a_{ij}=e^{\prime}a_{ji}$. For
$i=q_0, j=q+q_0-1$, we have $a_{ij}=1$ and $\omega^{i-j}=1$, and
the same holds for $i=q+q_0-1,j=q_0.$ These equations yield
$e=e^{\prime}=1.$ Hence $a_{ij}\omega^{q_0(i-j)}= a_{ij}=a_{ji}$.
Suppose now that $a_{ij}\neq 0$. Then
$\omega^{q_0(i-j)}=1$, and hence $i-j$ is divisible by $q-1$.
Since $i+j<q+2q_0-1$, this only leaves three possibilities, namely
$i=j$, $i=q-1+j$ and $j=q-1+i$. So,
$F(X^{\prime},Y^{\prime})=G_1(X^{\prime}Y^{\prime})+\ldots
a_{q-1+j,j}{X^{\prime}}^{q-1+j}{Y^{\prime}}^j+\ldots +
a_{i,q-1+i}{X^{\prime}}^i{Y^{\prime}}^{q-1+i}+ \ldots$. As
$a_{ij}=a_{ji}$, this gives
$$
\begin{array}{lcl}
F(X^{\prime},Y^{\prime}) & = & 
G_1(X^{\prime}Y^{\prime})+\ldots
a_{q-1+j,j}{X^{\prime}}^{q-1+j}Y^j+a_{j,q-1+j}{X^{\prime}}^j{Y^{\prime}}^{q-1+j}+ \\
{} & {} & +\ldots + a_{i,q-1+i}{X^{\prime}}^i{Y^{\prime}}^{q-1+i}+ 
a_{q-1+i,i}{X^{\prime}}^{q-1+i}Y^i+ \ldots\\
{} &  = & G_1(X^{\prime}Y^{\prime})+({X^{\prime}}^{q-1}+{Y^{\prime}}^{q-1})(a_{q-1+j,j}{X^{\prime}}^j{Y^{\prime}}^j
+a_{i,q-1+i}{X^{\prime}}^i{Y^{\prime}}^i+
\ldots)\\
{} & = &
G_1(X^{\prime}Y^{\prime})+G_2(X^{\prime}Y^{\prime})(X^{q-1}+Y^{q-1}).
\end{array}
$$
It remains to prove that $\deg(G_1)=\frac{1}{2}q-1$. Let $r$ be
the line of $\alpha$ with equation $X^{\prime}+Y^{\prime}=0$. From
(\ref{vPx+by}), the intersection divisor $\pi_{\ell}^{-1}(r)$ is
$(2q_0+1)P(1:1:0:1:1)+\sum
P(1:a:a^{q_0+1}+a^{q_0}:a:a^2+a^{2q_0}+a^{2q_0+2})$ with $a$
ranging over $\fq\setminus\{0,1\}$. Equivalently, for
$I=(\pi_{\ell}(\cX),r;\pi_\ell(P))$ we have
$$
 I = \left \{
\begin{array}{lll}
 1 & \mbox{\em for\ $P=(1:a:a^{q_0+1}+a^{q_0}:a:a^2+a^{2q_0}+a^{2q_0+2}); a\in\fq\setminus \{0,1\}$\/}, \\
 2q_0+1 & \mbox{\em for\ $P=(1:1:0:1:1)$\/}, \\
  0    & \mbox{\em otherwise.\/}
                         \end{array}
          \right.
$$
Since $\pi_\ell(1:1:0:1:1)$ is the infinite point of $r$ but
$\pi_{\ell}(1:a:a^{q_0+1}+a^{q_0}:a:a^2+a^{2q_0}:a^{2q_0+2})$ with
$a\in\fq\setminus\{0,1\}$ is a point at finite distance on $r$, we
deduce from the equation of $\pi_{\ell}(\cX)$ that
$\deg(G_1(T^2))=q-2$. More precisely, the roots of $G_1(T)$ are
the elements of the set $\{(a^q_0(a+1))^{-1}\}\mid a\in
\fq\setminus\{0,1\}$. To obtain such a set in a simpler form, we
note that $t\mapsto t^{2q_0-1}$ is a bijection of $\fq$. Putting
$u:=(a+1)^{-2q_0-1}$, we have
$(a^q_0(a+1))^{-1}=(u+u^{q_0})^{2q_0-1}.$ It turns out that the
roots of $G_1(T)$ can be also written in the form
$(u+u^{q_0})^{2q_0-1}$ with $u$ ranging over $\fq\setminus
\{0,1\}$.
     \end{proof}
To obtain an explicit expression for $G_1(T)$ and $G_2(T)$, we
need the following result from finite field theory, see
\cite[Section 1.4]{h}.
{\it An element $a\in \fq$ is of trace $0$ or trace $1$ according
as the polynomial $x^2+x+a$ is reducible or irreducible over
$\fq$.}
Let $Cl_0$ and $Cl_1$ be the set of all elements of trace $0$ and
of trace $1$ in $\fq$, respectively. Equivalently, $Cl_0$ and
$Cl_1$ consist of all roots of the polynomials
   $$
P_0(X):=X+X^2+X^4+\ldots+X^{q/2}\, ,\qquad
P_1(X)=1+X+X^2+X^4+\ldots+X^{q/2}\, ,
   $$
respectively. Furthermore,
$Cl_0$ is an additive subgroup of $\fq$ of index $2$ and its coset
is $Cl_1$. Also, it is easily seen that $Cl_0=\{c^{q_0}+c \ | \
c\in \fq\}$. Define
$$
G^-(T):=1+\sum_{i=0}^{s-1}T^{2^i(2q_0+1)-(q_0+1)}(1+T)^{2^i}\,, $$
and
$$ H(T):=TG^-(T)\,.
$$
\begin{lemma}\label{l17}
$H(T)$ is the polynomial with the lowest degree (i.e. $\deg
(H(T))=\frac{1}{2}q$) whose roots are the $(2q_0-1)$-{\rm th}
powers of the elements in $Cl_0$.
\end{lemma}
\begin{proof}
Since $Cl_0$ consists of $\frac{1}{2}q$ elements, and the map
$x\to x^{2q_0-1}$ in $\fq$ is bijective, we have to prove that if
$a\,\in\,Cl_0$, then $a^{2q_0-1}$ is a root of $H(T)$. Let $a\in
Cl_0$, and put $t=a^{2q_0-1}$. Then $t^q=t$, $t^{q_0+1}=a^{q_0}$,
$t^{2q_0+1}=a$, $t^{q/2+q_0}=a^{q/2}$. Now, $$
t^{q_0}H(t)=t^{q_0+1}+\sum_{i=0}^{s-1}t^{(2q_0+1)2^i}+
\sum_{i=0}^{s-1}t^{(2q_0+2)2^i}= $$
$$ =a^{q_0}+\sum_{i=0}^{s-1}a^{2^i}+\sum_{i=0}^{s-1}a^{(2q_0)2^i}=
$$
$$ =a^{q_0}+(a+a^2+\ldots+a^{q_0/2})+ (a^{2q_0}+a^{4q_0}+\ldots
+a^{q/2})=P_0(a)=0\,. $$
\end{proof}
\begin{lemma}\label{ciclicq-1D}
$G_1(T)=eG^-(T)$ with $e\in \fq$.
\end{lemma}
\begin{proof}
As we have seen in the final part of the proof of Proposition
\ref{propo-1}, $G_1(T)$ has degree $\frac{1}{2}q-1$ because its
roots are $(u+u^{q_0})^{2q_0-1}$ with $u$ ranging over
$\fq\backslash\{0,1\}$. In particular, every root of $G_1(T)$ is
obtained exactly twice, as
$(u+u^{q_0})^{2q_0-1}=(v+v^{q_0})^{2q_0-1}$ happens if and only if
either $v=u$ or $v=u+1$. On the other hand, such elements are
precisely the $2q_0-1$-powers of the non-zero elements in $Cl_0$.
By Lemma \ref{l17}, we obtain $TG_1(T)=eH(T))$ with $e\neq 0$,
whence the claim follows.
\end{proof}
\begin{lemma}
\label{cyclicq-1A} $G_2(T)=1+T^{q_0}$ and  $e=1$.
\end{lemma}
\begin{proof}
As we have already noted, $x$ is a local parameter at $A_0$, and
the local expansions of $by$ at $A_0$ is
$by=x^{2q_0+1}+x^{q+2q_0}+\ldots$. Thus, by Equation
(\ref{vxq0+1+yq0})
$$\xi\eta=\frac{x^{2q_0+2}+x^{q+2q_0+1}+\ldots}{x^{2q_0+2}+\ldots}=1+x^{q-1}+\ldots.$$
Write $G_2(T)=a_0+a_1T^m+\ldots$
where $m=2^nk$, $k$ odd, and $m\leq q_0$. Then
$G_2(\xi\eta)=a_0+a_1(1+x^{q-1}+\ldots)^m+\ldots=a_0+a_1+a_1x^{2^n(q-1)}+\ldots$,
where $\ldots$ indicate terms of degree greater than $2^n(q-1)$.
By Proposition \ref{propo-1} and Lemma \ref{ciclicq-1D},
$G_1(\xi\eta)+G_2(\xi\eta)(\xi^{q-1}+\eta^{q-1})=0$, and
$G_1(\xi\eta)=e(1+\ldots)$. Thus,
$v_{A_0}(G_2(\xi\eta)(\xi^{q-1}+\eta^{q-1})=
v_{A_0}((a_0+a_1)x^{-q_0(q-1)}+a_1x^{2^n(q-1)}x^{-q_0(q-1)}+\ldots)=
0$. This yields $a_0=a_1$ and $2^n=q_0$. Hence
$G_2(T)=a_0(1+T^{q_0})$. Moreover,
$G_1(\xi\eta)+G_2(\xi\eta)(\xi^{q-1}+\eta^{q-1})$ becomes
$e+a_0+\ldots$ where $\ldots$ indicate terms of positive degrees.
Hence $e=a_0$.
\end{proof}
As a corollary to Proposition \ref{propo-1} and to the previous
two lemmas, we have that an equation for $\pi_\ell(\cX)$ is
$$G^-(X^\prime Y^\prime)=
({X^\prime}^{q-1}+{Y^\prime}^{q-1})((X^\prime Y^\prime)^{q_0}
+1).$$ So we have found the desired plane model on which both $\bf
N$ and $[W]$ act linearly:
\begin{theorem}\label{thm4.6}
Let $G^-(T)\in \mathbb{F}_q[T]$ be defined to be
$$
G^-(T)=1+\sum_{i=0}^{s-1}T^{2^i(2q_0+1)-(q_0+1)}(1+T)^{2^i}.
$$
Then $\cX$ is birationally $\mathbb{F}_{q}$-isomorphic to the
plane curve of equation
\begin{equation}
\label{gminus} G^-(XY)= (X^{q-1}+Y^{q-1})((XY)^{q_0} +1).
\end{equation}
\end{theorem}
We are in a position to give an explicit equation for the quotient
curve $\cX_{\cU}$.
\begin{theorem}
\label{ff2} For every divisor $r$ of $q-1$, the quotient curve of
the DLS-curve associated to the cyclic subgroup $\cU$ of
${\cS}z(q)$ of order $r$ has genus $g=g_{\cU}=\frac{1}{r}q_0(q-1)$
and is $\fq$-isomorphic to the non-singular model of the plane
curve of equation
$$
Y^{(q-1)/r}\Big(1+\sum_{i=0}^{s-1}X^{2^i(2q_0+1)-(q_0+1)}(1+X)^{2^i}\Big)
=(X^{q_0}+1)(Y^{2(q-1)/r}+X^{q-1}).
$$
\end{theorem}
\begin{proof}
Let $\varphi_r:\pi_\ell(\cX) \mapsto
{\mathbf{P}}^2(\bar{\mathbf{F}}_q)$ be the rational map
$\varphi_r:=(1:X^{\prime}:Y^{\prime})\mapsto
(1:X^{\prime}Y^{\prime}:{Y^{\prime}}^r)$. Given a point
$Q:=(1:u:v)\in \Im (\varphi_r)$ with $v \neq 0$, let
$P:=(1:u':v')\in \varphi_r^{-1}(Q)$. For $i=1,\ldots,r$, let
$P_i:=(1:\tau^{-q_0i}u':\tau^{q_0i} v')$ with $\tau$ an element of
order $(q-1)/r$ in the multiplicative group of $\fq$. Then $P_i\in
\varphi^{-1}(Q)$. On the other hand, the equation
$v={Y^{\prime}}^r$ has precisely $r$ solutions, namely
$\tau^{q_0i} v'$ with $i=1,\ldots,r$. Hence $\varphi_{r}$ has
degree $r$, and
$\varphi_r^{-1}(Q)=\{(1:\tau^{-q_0i}u/v':\tau^{q_0i}
v')|i=1,\ldots,r\}$. This together with the third claim in Lemma
\ref{cyclicq-1} shows that a non--singular model of
$\varphi_r(\pi_\ell(\cX))$ is the quotient curve of $\cX$
associated to the subgroup ${\cU}$. Finally, Theorem \ref{thm4.6}
together with a direct computation gives the desired equation.
\end{proof}

\section{Quotient curves arising from the Singer type subgroups}\label{singertype}

{}From the classification of the subgroups of $\mathcal{S}z(q)$, see
Result \ref{Ris4} in Section \ref{prel}, there exist two cyclic
groups of Singer type up to conjugacy in
$\Aut(\cX)\cong\mathcal{S}z(q)$, one of order $q+2q_0+1$, the long
Singer subgroup $ \D^+$, and one of order $q-2q_0+1$, the short
Singer subgroup ${\D}^-$.
We look for two (singular) plane curves $\cD^+$ and $\cD^-$, both
birationally isomorphic to $\cX$, such that ${\D}^+$ acts on
$\cD^+$ and ${\D}^-$ acts on $\cD^-$ as a linear collineation
group. We use the same approach as in Section \ref{ordq-1}. For
this purpose, we choose a generator ${\bf g}$ of ${ \D}^+$ (or
${\D}^-$) represented by a $4\times4$-matrix $B$ over
$\mathbb{F}_q$, and check that  ${\bf g}$ has exactly five fixed
points: one, say $B_0$, is defined over $\mathbb{F}_q$ and the
other four, say $B_1,\ldots,B_4$ are defined over
$\mathbb{F}_{q^4}$. We can arrange the indices in such a way that
$B_2={\bf F}(B_1), \ B_3={\bf F}^2(B_1), \ B_4={\bf F}^3(B_1)$,
where ${\bf F}$ denotes the semilinear collineation
$(X_0:X_1:X_2:X_3:X_4)\mapsto (X_0^q:X_1^q:X_2^q:X_3^q:X_4^q)$. The
point $B_0$ is not on $\cX$, while $B_1$ (and hence each of
$B_2,B_3,B_4$) belongs to $\cX$ or does not according as ${\bf g}$
generates $\D^+$ or $\D^-$. In both cases, the points $B_i$ are
linearly independent. Every element of order four in the
normaliser of the Singer group generated by ${\bf g}$, fixes $B_0$
and preserves the set $\{B_1,B_2,B_3,B_4\}$. More precisely, it
acts on $\{B_1,B_2,B_3,B_4\}$ as either $(B_1B_2B_3B_4)$, or
$(B_1B_4B_3B_2)$. Each element of order two in the normaliser
changes $B_1$ with $B_3$, and $B_2$ with $B_4$. At this point, we
note that the line $\ell$ through $B_2$ and $B_4$ is defined over
$\mathbb{F}_{q^2}$ as it is left invariant by ${\bf F}^2$. The
linear system of all hyperplanes through $\ell$ cuts out on $\cX$
a two-dimensional linear series $g_2^n$ defined over
$\mathbb{F}_{q^2}$. The degree $n$ is either $q+2q_0+1$ or
$q+2q_0-1$, according as ${\bf g}$ generates $\D^+$ or $\D^-$. The
irreducible plane curve $\cD$ associated with $g_2^n$ is left
invariant by the Singer group generated by ${\bf g}$ which acts on
it as a linear collineation group. It turns out that $\cD^+=\cD$
for the long Singer subgroup $\bf D^+$ and $\cD^-=\cD$ for the
short Singer subgroup $\bf D^-$ have the desired properties.
To find an explicit equation for $\cD$ we transform the above
matrix $B$ into its diagonal form $\Lambda$ defined over
$\mathbb{F}_{q^4}$. The eigenvalues are $1$ and $\lambda^{q^i}$,
$i=0,1,2,3$ for an element $\lambda \in{\mathbb{F}}_{q^4}$ whose
order is either $q+2q_0+1$ or $q-2q_0+1$ according as ${\bf g}$
generates $\D^+$ or $\D^-$. Once we have chosen a new frame
$(X_0^{\prime}:X_1^{\prime}:X_2^{\prime}:X_3^{\prime}:X_4^{\prime})$
whose fundamental simplex is $B_0B_1B_2B_3B_4$, the plane curve
$\cD$ turns out to be the projection $\pi_\ell$ of $\cX$ from the
line $\ell$ viewed as the vertex, to the plane through
$B_0B_1B_3$. In other terms, the projection is
$\pi_{\ell}: (X_0^{\prime}:X_1^{\prime}:X_2^{\prime}:X_3^{\prime}:X_4^{\prime})
\mapsto (X_0^{\prime}:X_1^{\prime}:X_3^{\prime})$. The equation of
$\cD$ is given in the following theorems.
\begin{theorem}
\label{modellong} Let $G^+(T)\in \mathbb{F}_q[T]$ be defined to be
$$ G^+(T)=1+\sum_{i=0}^{s-1}
T^{2^iq_0}(1+T)^{2^i(q_0+1)-q_0}+T^{q/2}\,. $$ Then $\cX$ is
birationally $\mathbb{F}_{q^2}$-isomorphic to the plane curve
$\cD^+$ of equation $$G^+(XY)=X^{q+2q_0+1}+Y^{q+2q_0+1}.$$
\end{theorem}
\begin{theorem}
\label{modelsmall} Let $G^-(T)\in \mathbb{F}_q[T]$ be defined to
be $$
G^-(T)=1+\sum_{i=0}^{s-1}T^{2^i(2q_0+1)-(q_0+1)}(1+T)^{2^i}\,. $$
Let
$b:=\lambda^{q_0}+\lambda^{q_0-1}+\lambda^{-q_0}+\lambda^{-(q_0-1)}$
for an element $\lambda \in \mathbb{F}_q^4$ of order $q-2q_0+1$.
Then $b\in\mathbb{F}_q$, and $\cX$ is birationally
$\mathbb{F}_{q^2}$-isomorphic to the plane curve $\cD^-$ of
equation $$bG^-(XY)= (X^{q-2q_0+1}+Y^{q-2q_0+1})((XY)^{q_0-1}
+(XY)^{2q_0-1}).$$
\end{theorem}
In carrying out the necessary computations for the proof, we will
need to use some more notation. Fix an element $\lambda$ in the
multiplicative group of $\mathbb{F}_{q^4}^{\star}$ whose order is
either $q+2q_0+1$ or $q-2q_0+1$ ; that is either
$\lambda=w^{(q^2-1)(q-2q_0+1)}$, or
$\lambda=w^{(q^2-1)(q+2q_0+1)}$, where $w$ is a primitive element
of $\mathbb{F}_{q^4}$. Also, let
$$
b:=\begin{cases}
\lambda^{q_0}+\lambda^{q_0+1}+\lambda^{-q_0}+\lambda^{-(q_0+1)}  &
\text{for $\lambda^{q+2q_0+1}=1,$}\\
\lambda^{q_0}+\lambda^{q_0-1}+\lambda^{-q_0}+\lambda^{-(q_0-1)}
                      & \text{for $\lambda^{q-2q_0+1}=1$.}
\end{cases}
$$
$$
 \mu:=\frac{(1+\lambda)^2}{b\lambda}.
$$
$$
\rho:=\frac{1}{ b^{q_0} \mu}.
$$
Furthermore, for $x,y\in \fq(\cX)$ as in Proposition \ref{functionfield} let
$$
h(x,y):=bxy+x^{2q_0+2}+b^{2q_0}y^{2q_0}\,,
$$
and
$$
k(x,y):=b^{q_0-1}x+b^{q_0}y+x^{q_0+1}+b^{q_0}y^{q_0}\,.
$$
The following equalities are straightforward to check.
$$
b^{2q_0}=(\lambda+\lambda^{-1})^q+(\lambda+\lambda^{-1})\,,
$$
$$
\lambda^4+b^{2q_0}\lambda^3+b^2\lambda^2+b^{2q_0}\lambda+1=0\,,
$$
$$
\lambda^2=1+\lambda\mu b\,,
$$
$$
\mu^{q+1}=1, {\rm \,\, and \,\, hence \,\, }
\mu\in\mathbb{F}_{q^2},
$$
$$
h(x,y)^{q_0}=b^{q_0}x^{q_0}y^{q_0}+x^{2q_0+1}+by\,,
$$
$$
h(x,y)^q=bx^qy^q+x^{2q+2q_0}+b^{2q_0}y^{2q_0}=x^{q+2q_0+1}+bx^qy+b^{2q_0}y^{2q_0}\,,
$$
$$
(x^{q_0+1}+b^{q_0}y^{q_0})h(x,y)^{q_0}+b^{q_0+1}y^{q_0+1}+x^{q_0}h(x,y)=0\,,
$$
$$
k(x,y+1)=k(x,y)\,,
$$
\begin{equation}
\label{3.4)}h\Big(\frac{by}{h(x,y)},\frac{x}{bh(x,y)}\Big)=h(x,y)^{-1}\,,
\end{equation}
\begin{equation}
\label{3.6)}k\Big(\frac{by}{h(x,y)},\frac{x}{bh(x,y)}\Big)=k(x,y)/h(x,y)\,.
\end{equation}
We will also need a technical lemma.
Let
$P:=(1:u:u^{q_0+1}+b^{q_0}y^{q_0}:bv:buv+u^{2q_0+2}+b^{2q_0}y^{2q_0})$
be a point of $\cX$. Then $t:=x+u$ is a local parameter at $P$.
\begin{lemma}
There exists ${\bar y}\,\in\,{\bfq}(\cX)$ such that
$v_{P} (\bar y)\ge q$ and
\begin{equation}\label{ramo}
by=bv+u^{2q_0}t+ (u+v^q)t^{2q_0}+t^{2q_0+1}+{\bar y}\,.
\end{equation}
\end{lemma}
\begin{proof} From $x^{2q_0}(x^q+x)=by^q+by$ we have $$by=by^q+(t+u)^{2q_0}(t^q+
t+u^q+u)\,, $$ whence the claim follows for ${\bar y}:=b(y-v)^q$.
\end{proof}
Let $B:=WT_{o,b}$, that is $$ B:=
 \left( \begin{array}{ccccc}
0                      & 0            & 0 & 0 & 1 \\ 0
& 0            & 0 & 1 & 0 \\ 0                      & 0
& 1 & 0 & b^{q_0} \\ 0                        & 1            & 0 &
0 & b \\
                      1& 0 & 0 & b & b^{2q_0}
\end{array} \right)\,.
$$
Let $$ M:=
 \left( \begin{array}{ccccc}
0 & b^{q_0-1} & 1 & b^{q_0-1} & 0\\ \mu & 1 & 0 & \lambda &
\lambda\mu\\ \mu^q & 1 & 0 & \lambda^q & \lambda^q\mu^q\\
\mu^{q^2} & 1 & 0 & \lambda^{q^2} & \lambda^{q^2}\mu^{q^2}\\
\mu^{q^3} & 1 & 0 & \lambda^{q^3} & \lambda^{q^3}\mu^{q^3}
\end{array} \right)\,,\quad
{\rm and} \quad \Lambda:=
 \left( \begin{array}{ccccc}
 1 & 0 & 0 & 0 & 0\\
 0 & \lambda & 0 & 0& 0\\
 0 & 0 & \lambda^{q}& 0 & 0\\
 0 & 0& 0 & \lambda^{-1}& 0\\
 0& 0 & 0 & 0 & \lambda^{-q}
\end{array} \right) \,.
$$ A straightforward computation shows that
\begin{equation}\label{lem5.4}
MBM^{-1}=\Lambda.
\end{equation}
Let $[M]$ be the collineation associated with $M$. Then $[M]$ is a
morphism $\cX \mapsto \P^4({\mathbb{F}}_{q^4})$, and $[\Lambda]$
is a linear automorphism of $\cY:=[M](\cX)$. The algebraic curve
$\cY$ can be viewed as a parameterised curve
associated to the morphism $g$ with coordinate functions
  $$
g:= (g_0:g_1:g_2:g_3:g_4)
  $$
where
  $$
g_0:=b^{q_0-1}f_1+f_2+b^{q_0-1}f_3=k(x,y)\,, $$ $$ g_1:=\mu
f_0+f_1+\lambda f_3+\lambda \mu f_4=\mu+x+\lambda by + \mu\lambda
h(x,y)\,, $$ $$ g_2:=\mu^{q} f_0+f_1+\lambda^q f_3+\lambda^q
\mu^{q} f_4=\mu^{-1}+x+\lambda^q by +
 \mu^{-1}\lambda^q h(x,y)\,,
$$ $$ g_3:= \mu^{q^2} f_0+f_1+\lambda^{q^2} f_3+\lambda^{q^2}
\mu^{q^2} f_4= \mu+x+\lambda^{-1} by + \mu\lambda^{-1} h(x,y)\,,
$$ $$ g_4:=\mu^{q^3} f_0+f_1+\lambda^{q^3} f_3+\lambda^{q^3}
\mu^{q^3} f_4= \mu^{-1}+x+\lambda^{-q} by + \mu^{-1}\lambda^{-q}
h(x,y) \,. $$
Note that the fixed points of $[\Lambda]$ are $B_0:=(1:0:0:0:0)$,
$B_1:=(0:1:0:0:0)$, $B_2:=(0:0:1:0:0)$, $B_3:=(0:0:0:1:0)$ and
$B_4:=(0:0:0:0:1)$.
The following lemma follows from Equation (\ref{ramo}) together
with a straightforward computation.
\begin{lemma}\label{expansion}
For a point
$P:=(1:u:u^{q_0+1}+b^{q_0}y^{q_0}:bv:buv+u^{2q_0+2}+b^{2q_0}y^{2q_0})\in\cX$
the local expansion of the coordinate functions of $g$ at $P$ are
\[
\begin{array}{ll}
g_0= & g_0(P)+t(b^{q_0-1}(1+u^{2q_0}) +u^{q_0})+t^{q_0}(u+u^{q})\\
{} &+t^{q_0+1}+t^{2q_0}(b^{q_0-1}(u+u^q))+t^{2q_0+1}b^{q_0-1}
+{\bar g_0}\,,
\end{array}
\]
\[
\begin{array}{ll}
g_1= & g_1(P)+t(1+\lambda (u^{2q_0}+\mu u^{2q_0+1} +\mu bv))\\
{} & +t^{2q_0}\lambda(u+u^q)(1+\mu u^q)+t^{2q_0+1}\lambda (1+\mu
u^q)+{\bar g_1}\,,
\end{array}
\]
\[
\begin{array}{ll}
g_2= & g_2(P)+t(1+\lambda^q (u^{2q_0}+\mu^{-1} u^{2q_0+1}
+\mu^{-1} bv))\\
{} & +t^{2q_0}\lambda^q (u+u^q)(1+\mu u^q) +t^{2q_0+1}\lambda^q
(1+\mu^{-1} u^q) +{\bar g_2}\,,
\end{array}
\]
\[
\begin{array}{ll}
g_3= & g_3(P)+t(1+\lambda^{-1} (u^{2q_0}+\mu u^{2q_0+1} +\mu
bv))\\
{} & +t^{2q_0}\lambda^{-1} (u+u^q)(1+\mu u^q)
 +t^{2q_0+1}\lambda^{-1} (1+\mu u^q)
+{\bar g_3}\,,
\end{array}
\]
\[
\begin{array}{ll}
g_4= & g_4(P)+t(1+\lambda^{-q} (u^{2q_0}+\mu u^{2q_0+1} +\mu
bv))\\
{} & +t^{2q_0}\lambda^{-q} (u+u^q)(1+\mu u^q)
+t^{2q_0+1}\lambda^{-q} (1+\mu u^q) +{\bar g_4}\,,
\end{array}
\]
with $v_P({\bar g_i})\ge q$ for $i=0,\ldots, 4$.
\end{lemma}
We will also use a result on finite fields.
\begin{lemma}\label{l0}
The system in $T$
\begin{itemize}
\label{sys}
\item [1)]  $b^{q_0}T^{q_0}+b^{q_0}T=b^{q_0-1}\mu+\mu^{q_0+1},$
\item [2)]  $b(T^q+T)=\mu^{2q_0}(\mu^q+\mu)$
\end{itemize}
is not solvable in ${\bar{\mathbb{F}}}_q$ for
$\lambda^{q+2q_0+1}=1$, but it has exactly two solutions for
$\lambda^{q-2q_0+1}=1$, namely $(\mu/b) \lambda^q$ and $
(\mu/b)\lambda^q+1$.
\end{lemma}
\begin{proof}
Assume that the above system is consistent, and let $z$ denote a
solution. We show that
\begin{itemize}
\item [3)] $z^2+z=(\mu/b)^2.$
\end{itemize}
{}From 1) it follows
$bz^q+bz^{2q_0}+b^{1-2q_0}\mu^{2q_0}+\mu^{q+2q_0}=0$. This
together with 2) yield
$b^{2q_0}(z^{2q_0}+z)+\mu^{2q_0}+b^{2q_0-1}\mu^{2q_0+1}=0.$ Adding
it to the squared of 1)  gives
$b^{2q_0}z^2+b^{2q_0}z=\mu^{2q_0+1}(\mu+1/\mu+b^{2q_0-1})+b^{2q_0-2}\mu^2=0$.
Hence 3) follows form (1.1). Claim 3) implies that the system has
at most two solutions, $z$ and $z+1$. Actually, one of them is
$(\mu/b) \lambda^q$. In fact, $(\mu/ b)^2 \lambda^{2q}+(\mu /b)
\lambda^q+(\mu/b)^2=0$ holds if and only if $
\lambda^{2q}+\lambda^{q}b^2\lambda/(1+\lambda^2) = 1$. This is
true for $\lambda^{q^2+1}=1$ and (1.1). It is straightforward to
check that $(\mu/b) \lambda^q$ satisfies 1) and 2) if and only if
$\lambda^{q-2q_0+1}=1$.
\end{proof}
\begin{lemma}\label{int1}
For $\lambda^{q+2q_0+1}=1$, we have $B_1,B_3\not\in\cY$, and
$\ell\cap \cY=\emptyset$.
\end{lemma}
\begin{proof}
Let
$P:=(1:u:u^{q_0+1}+b^{q_0}y^{q_0}:bv:buv+u^{2q_0+2}+b^{2q_0}y^{2q_0})$
be a point of $\cX$. If $g(P)\in \ell$, then the equation
$\mu+u+\lambda bv + \mu\lambda
buv+u^{2q_0+2}+(bv)^{2q_0}=\mu+u+\lambda^{-1} bv + \mu\lambda^{-1}
buv+u^{2q_0+2}+(bv)^{2q_0}=0$ yields $bv=\mu
buv+u^{2q_0+2}+(bv)^{2q_0}$ and hence $u=\mu$. Now, suppose that
$g(P)=B_i$, $i=1,3$. From $\mu^{-1}+u+\lambda^q v +
\mu^{-1}\lambda^q
buv+u^{2q_0+2}+(bv)^{2q_0}=\mu^{-1}+u+\lambda^{-q} v +
\mu^{-1}\lambda^{-q} buv+u^{2q_0+2}+(bv)^{2q_0}=0$ we have
$v=\mu^{-1}buv+u^{2q_0+2}+(bv)^{2q_0}$ and $u=\mu^{-1}$. At this
point, it is enough to show that if $u\, \in \{\mu,\mu^{-1}\}$,
then $(u,v)$ is not a point of the plane curve $\cC^{\prime}$ of
equation $k(X,Y)=b^{q_0-1}X+b^{q_0}Y+X^{q_0+1}+(bY)^{q_0}=0$.
Actually, this claim follows from Lemma \ref{l0} as both $\cC_b$
and $\cC^{\prime}$ are defined over $\fq$.
\end{proof}
Let
$$P_1:=(1:1/\mu:(1/\mu)^{q_0+1}+(1/(\lambda\mu)+b)^{q_0}:1/(\lambda\mu)+b:
$$
$$\quad\qquad(1/\mu)(1/(\lambda\mu)+b) +(1/\mu)^{2q_0+2}+
(1/(\lambda\mu)+b)^{2q_0}) ,$$
$$P_2:=(1:\mu:(\mu)^{2q_0+1}+(\mu\lambda^q)^{q_0}:\mu\lambda^q:\mu^2\lambda^q+(\mu)^{2q_0+3}
\lambda^q+(\mu\lambda^q)^{2q_0}),$$
$$P_3:=(1:1/\mu:(1/\mu)^{q_0+1}+1/(\lambda\mu)^{q_0}:1/(\lambda\mu):1/(\lambda\mu^2)+(1/\mu)^{2q_0+2}
+1/(\lambda\mu)^{2q_0}),$$
$$P_4:=(1:\mu:\mu^{q_0+1}+(\mu\lambda^q+b)^{q_0}:\mu\lambda^q+b:\mu(\mu\lambda^q+b)+\mu^{2q_0+2}+
(\mu\lambda^q+b)^{2q_0}).$$ It is straightforward to check that
$P_i\, \in \, \cX$ for $i=1,\ldots,4$.
\begin{lemma}\label{int2}
For $\lambda^{q-2q_0+1}=1$, we have $g(P_i)=B_i$ for $i=1,
\ldots,4$. Furthermore, $B_2$ and $B_4$ are the common points of
$\ell$ and $\cY$.
\end{lemma}
\begin{proof}
A direct computation proves that $g(P_i)=B_i$, for $i=1,\ldots,4$.
Given a point
$P:=(1:u:u^{q_0+1}+(bv)^{q_0}:bv:ubv+x^{2q_0+2}+(by)^{2q_0})\in
\cX$, suppose that $g(P)\in \ell$. Then $\mu+u+\lambda bv +
\mu\lambda buv+u^{2q_0+2}+(bv)^{2q_0}=\mu+u+\lambda^{-1} bv +
\mu\lambda^{-1} buv+u^{2q_0+2}+(bv)^{2q_0}=0$ yields $bv=\mu
buv+u^{2q_0+2}+(bv)^{2q_0}$ and hence $u=\mu$. By Lemma \ref{l0}
$P\in\{P_2,P_4\}$, whence $g(P)\in \{B_2,B_4\}$.
\end{proof}
\begin{lemma}\label{preimage}
For $\lambda^{q-2q_0+1}=1$, we have
$$ v_{P_{1}}(g_0)=q_0, \qquad
v_{P_1}(g_1)= 0, \qquad v_{P_1}(g_3)=2q_0\,, $$ $$
v_{P_{2}}(g_0)=q_0, \qquad v_{P_2}(g_1)= q, \qquad
v_{P_2}(g_3)=1\,, $$ $$ v_{P_{3}}(g_0)=q_0, \qquad v_{P_3}(g_1)=
2q_0, \qquad v_{P_3}(g_3)=0\,, $$ $$ v_{P_{4}}(g_0)=q_0, \qquad
v_{P_4}(g_1)=1, \qquad v_{P_4}(g_3)= q\,. $$
\end{lemma}
\begin{proof}
Since the order-sequence of $\cX$ at any point $Q\in
\cX\setminus\cX(\fq)$ is $(0,1,q_0,2q_0,q)$, we have
$v_{P_i}(g_j)\le q$ for $i=1,\ldots, 4$ and $j=0,1,3$. Then the
lemma follows from Lemma \ref{expansion} together with some
computation.
\end{proof}
\begin{corollary}\label{image}
For $\lambda^{q-2q_0+1}=1$, $\pi_\ell(B_2)=(0:0:1)$ and
$\pi_\ell(B_4)=(0:1:0)$.
\end{corollary}
\begin{lemma}\label{divisor+}
For $\lambda^{q+2q_0+1}=1$, we have
\begin{equation}
 v_{P}(g_0) \ = \left \{ \begin{array}{lll}
 1 & \mbox{\em for $P=[B]^jA_0$ with
 $j=0,\ldots, q+2q_0, \, j\neq 1$\/}, \\
-q-2q_0 & \mbox{\em for $P=A_4$\/}, \\
 0    & \mbox{\em otherwise.}
                         \end{array}
          \right.
\end{equation}
\end{lemma}
\begin{proof}
The claim for $A_0$ and $A_4$ follows from Equations (\ref{vPx}),
(\ref{vPy}) and (\ref{vxq0+1+yq0}). Note that $A_4$ is the only
pole of $g_0$, hence $\ord (g_0)=q+2q_0$. By Equation
(\ref{lem5.4}), for any integer $j$ the point $[B]^jA_0$ is a zero
of $g_0$ apart from the case where $[B]^jA_0= A_4$. On the other
hand, $[B]^jA_0=A_4$ only holds for $j\equiv 1 \pmod{ q+2q_0+1}$,
and this completes the proof.
\end{proof}
\begin{lemma}\label{divisor-}
For $\lambda^{q-2q_0+1}=1$, we have
\begin{equation}
 v_{P}(g_0) \ = \left \{ \begin{array}{lll}
 q_0 & \mbox{\em for $P\in \{P_1,P_2,P_3,P_4\}$} \\
 1 & \mbox{\em for $P=[B]^jA_0$ with
 $j=0,\ldots, q-2q_0, \, j\neq 1$\/}, \\
-q-2q_0 & \mbox{\em for $P=P_\infty$\/}, \\
 0    & \mbox{\em otherwise.}
                         \end{array}
          \right.
\end{equation}
\end{lemma}
\begin{proof}
The assertion for $P_i$, $i=1, \ldots, 4$ follows from Lemma
\ref{preimage}. For the remaining cases, the proof is similar to
that of Lemma \ref{divisor+}.
\end{proof}
For $P\in \cX$, let $e'_P:=-\min\{v_P(g_0),v_P(g_1),v_P(g_3)\}$.
\begin{lemma}\label{E+}
For $\lambda^{q+2q_0+1}=1$ we have
\begin{equation}
 e'_P \ = \left \{ \begin{array}{lll}
 q+2q_0+1 & \mbox{\em for $P=A_4$} \\
 0    & \mbox{\em otherwise.}
                         \end{array}
          \right.
\end{equation}
\end{lemma}
\begin{proof}
$e^{\prime}_{A_4}=q+2q_0+1$ can be easily checked once the
valuations of the functions $f_i$, $i:=0,\ldots,4$ at
$A_4=P_{\infty}$ are computed. Note that, because of the setup of
the present paper, such computations are done in Section
\ref{nontame}. For $P\neq A_4$, the claim follows from Lemma
\ref{int1}.
\end{proof}
\begin{lemma}\label{E-}
For $\lambda^{q-2q_0+1}=1$ we have
\begin{equation}
 e'_P \ = \left \{ \begin{array}{lll}
 q+2q_0+1 & \mbox{\em for $P=A_4$} \\
 -1 & \mbox{\em for $P\in \{P_2,P_4\}$} \\
 0    & \mbox{\em otherwise.}
                         \end{array}
          \right.
\end{equation}
\end{lemma}
\begin{proof}
Similar to the previous one. For $P\neq A_4$, the claim follows
from Lemmas \ref{int2} and \ref{preimage}.
\end{proof}
The homogenous coordinates
$(X^{\prime}_0:X^{\prime}_1:X^{\prime}_3)$ provide a natural
projective frame in $\alpha$ with fundamental triangle
$B_0B_1B_3$. To write an equation of $\cD=\pi_\ell (\cY)$ we will
use in $\alpha$ the affine frame $(X^{\prime},Y^{\prime})$ arising
from the above projective frame by $X^{\prime}=\rho
\frac{X^{\prime}_1}{X^{\prime}_0}$, $Y^{\prime}=\rho
\frac{X^{\prime}_3}{X^{\prime}_0}$. Then $\cD$ has equation
$F(X^\prime,Y^\prime)=0$ where $F(X^\prime,Y^\prime)$ is an
absolutely irreducible polynomial with coefficients in $\fq$ which
satisfies $F(\xi,\eta)=0$, where $\xi,\eta$ are defined to be $\xi
:=\rho\,\frac{\mu + x+ \lambda by +  \lambda \mu h(x,y)
}{k(x,y)}$, $\eta := \rho\,\frac{\mu + x+ \lambda^{-1} by  +
\lambda^{-1} \mu h(x,y) }{k(x,y)}$.
\begin{lemma}
\label{degrees} The plane curve $\cD$ is birationally ${\mathbb
F}_{q^4}$-isomorphic to $\cX$. The degree of $\cD$ is $q+2q_0+1$
or $q+2q_0-1$ according as $\lambda^{q+2q_0+1}=1$ or
$\lambda^{q-2q_0+1}=1$.
\end{lemma}
\begin{proof}
We argue as in the proof of Lemma \ref{cyclicq-1}. Note that
$\pi_\ell g (A_0)=(0:1:1)$. Take any point
$Q:=(1:u:u^{q_0+1}+b^{q_0}v^{q_0}:bv:buv+u^{2q_0+2}+b^{2q_0}v^{2q_0})$
on $\cX$ such that $\pi_\ell g(Q)=\pi_\ell g(P)$. Then $g(Q)$ lies
on the hyperplane $X_0^\prime=0$, and $g_1(Q)=g_3(Q)$. By Lemmas
\ref{divisor+} and \ref{divisor-} together with straightforward
computation it turns out that $P=Q$. Moreover, by Lemma
\ref{expansion} for $u=v=0$, $1/\xi$ is a local parameter at $P$. From these facts we deduce that $\pi_{\ell}$ is birational.
Finally, Lemma \ref{E+} for $\lambda^{q+2q_0+1}=1$ and Lemma
\ref{E-} for $\lambda^{q-2q_0+1}=1$ imply the assertion concerning
the degree of $\cD$.
\end{proof}
The linear transformations $[\Lambda]$, $[MT_{o,b}M^{-1}]$ and
$[MWM^{-1}]$ preserve the line $\ell$. Hence, they act on the
set of planes through $\ell$, and give rise to linear
automorphisms of $\cD$. More precisely, the following lemmas hold.
\begin{lemma}\label{l2}
The automorphism $[T_{0,b}]$ acts on $\cD$ as the linear
transformation $(\xi,\eta)\mapsto(\lambda^2\eta,\lambda^{-2}\xi).$
\end{lemma}
\begin{proof}
The lemma follows from the following two relations
$$\frac{x+\mu+\lambda(by+b)+\lambda\mu h(x,y+1)}{k(x,y+1)}=
\lambda^2 \ \frac{x+\mu+\lambda^{-1}by+\lambda^{-1}\mu
h(x,y)}{k(x,y)}\,,$$
$$\frac{x+\mu+\lambda^{-1}(by+b)+\lambda^{-1}\mu
h(x,y+1)}{k(x,y+1)}= \lambda^{-2}\  \frac{x+\mu+\lambda
by+\lambda\mu h(x,y)}{k(x,y)}.$$
\end{proof}
\begin{lemma}\label{l3}
The automorphism $[W]$ acts on $\cD$ as the linear transformation
$(\xi,\eta)\mapsto(\lambda \eta,\lambda^{-1}\xi).$
\end{lemma}
\begin{proof}
The lemma follows from the following two relations, which are a
consequence of  Equations (\ref{3.4)}) and (\ref{3.6)}):
$$ \frac{\mu+\frac{by}{h(x,y)}+\lambda(\frac{x}{h(x,y)}+\mu
h(\frac{by}{h(x,y)},\frac{b^{-1}x}{h(x,y)}))}{k(\frac{by}{h(x,y)},
\frac{b^{-1}x}{h(x,y)})}= \lambda \ \frac{x+\mu+\lambda^{-1}(y+\mu
h(x,y))}{k(x,y)}, $$
$$ \frac{\mu+\frac{by}{h(x,y)}+\lambda^{-1}(\frac{x}{h(x,y)}+\mu
h(\frac{by}{h(x,y)},\frac{b^{-1}x}{h(x,y)}))}{k(\frac{by}{h(x,y)},
\frac{b^{-1}x}{h(x,y)})}= \lambda^{-1} \ \frac{x+\mu+\lambda(y+\mu
h(x,y))}{k(x,y)}\,. $$
\end{proof}
The following corollaries are straightforward to check.
\begin{corollary}\label{c1}
The automorphism $[B]$ acts on $\cD$ as the linear transformation
$(\xi,\eta)\mapsto (\lambda^{-1}\xi,\lambda\eta)$. In particular,
such an automorphism has order either  $q+2q_0+1$ or $q-2q_0+1$
according as $\lambda^{q+2q_0+1}=1$ or $\lambda^{q-2q_0+1}=1$.
\end{corollary}
\begin{corollary}\label{c2}
The automorphism $[BW]$ acts on $\cD$ as the linear transformation
$(\xi,\eta)\mapsto (\eta,\xi)$.
\end{corollary}
Corollary \ref{c1} is the essential tool for the proof of the
following result.
\begin{proposition}\label{propo0}
For $\lambda^{q+2q_0+1}=1$, the equation of $\cD$ can be written
in the form $$
G(X^{\prime}Y^{\prime})={X^{\prime}}^{q+2q_0+1}+{Y^{\prime}}^{q+2q_0+1},
$$ with $G(T)\in {\mathbb F}_{q^4}[T]$.
\end{proposition}
\begin{proof} We argue as in the proof of Proposition
\ref{propo-1}. Write the equation $F(X^{\prime},Y^{\prime})$ of
$\cD$ as $\sum a_{ij}X^{\prime}Y^{\prime}=0$. By Lemma
\ref{degrees} $F(X^{\prime},Y^{\prime})$ has degree $q+2q_0+1$.
Let $r_\infty$ be the line of $\alpha$ of equation $X_0^\prime=0$.
By Lemmas \ref{divisor+} and \ref{E+}, the intersection divisor of
$r_\infty $ is $(\pi_\ell g)^{-1}(r_\infty)=\sum_{j=0}^{q+2q_0}
[B]^jA_0=\sum_{j=0}^{q+2q_0} [M]^{-1}(0:\lambda^j \mu:\lambda^{qj}
\mu^{-1}: \lambda^{-j}\mu : \lambda^{-qj}\mu^{-1} )$. Then the
intersection between $\cD$ and $r_{\infty}$ consists of the points
$(0:\lambda^{j}:\lambda^{-j})$ where $j=0, \ldots, q+2q_0$. Taking
$\lambda^{q+2q_0+1}=1$ into account, this yields that the
monomials of degree $q+2q_0+1$ of $F(X^{\prime},Y^{\prime})$ are
${X^{\prime}}^{q+2q_0+1}$ and ${Y^{\prime}}^{q+2q_0+1}$. Since by
Corollary \ref{c1} the linear transformation $(\xi,\eta)\mapsto
(\lambda^{-1}\xi,\lambda\eta)$ fixes $\cD$, there exists $e\in
{\bar \fq}$ such that $$ \sum
a_{ij}{X^{\prime}}^i{Y^{\prime}}^j=e\sum a_{ij}\lambda^{i-j}
{X^{\prime}}^i{Y^{\prime}}^j\,, $$ that is
$a_{ij}=ea_{ij}\lambda^{i-j}$. Letting $i=q+2q_0+1,j=0$  yields
$e=1$. Furthermore, $a_{ij}\neq 0$ yields $\lambda^{i-j}=1$, that
is $i-j$ is divisible by $q+2q_0+1$. Since $i+j\leq q+2q_0+1$,
this only leaves three cases, namely $i=j$; $i=0,j=q+2q_0+1$;
$i=q+2q_0+1,j=0$.
\end{proof}
\begin{proposition}\label{propo1prelim} For $\lambda^{q-2q_0+1}=1$, the equation of $\cD$ can be
written in the form $$
G_1(X^{\prime}Y^{\prime})=G_2(X^{\prime}Y^{\prime})
({X^{\prime}}^{q-2q_0+1}+{Y^{\prime}}^{q-2q_0+1})\,, $$ with
$G_1(T),G_2(T)\in {\mathbb F}_{q^4}[T]$, and $\deg (G_2)=2q_0-1$.
\end{proposition}
\begin{proof}
Using the same notation as in the previous proof, from Lemmas
\ref{divisor-} and \ref{E-} the intersection divisor of $r_\infty
$ is $$(\pi_\ell
g)^{-1}(r_\infty)=q_0(P_1+P_3)+(q_0-1)(P_2+P_4)+\sum_{j=0}^{q-2q_0}
[B]^jA_0.$$ Then the intersection between $\cD$ and $r_{\infty}$
consists of $(0:0:1)$ and $(0:1:0)$, both counted $2q_0-1$ times,
together with the points $\{(0:\lambda^{j}:\lambda^{-j})\mid j=0,
\ldots, q-2q_0\}$, each counted just once. Since $j^{q-2q_0+1}=1$,
this implies that the monomials of degree $q+2q_0-1$ of
$F(X^{\prime},Y^{\prime})$ are
${X^{\prime}}^{q-2q_0+1}(X^{\prime}Y^{\prime})^{2q_0-1}+
{Y^{\prime}}^{q-2q_0+1}(X^{\prime}Y^{\prime})^{2q_0-1}$. By
Corollaries \ref{c1} and \ref{c2} both linear transformations
$(\xi,\eta)\mapsto (\lambda^{-1}\xi,\lambda\eta)$ and
$(\xi,\eta)\mapsto (\eta,\xi)$ fix $\cD$. Hence there exist $e\, ,
e'\, \in {\bar \fq}$ such that $$ \sum
a_{ij}{X^{\prime}}^i{Y^{\prime}}^j=e\sum a_{ij}\lambda^{i-j}
{X^{\prime}}^i{Y^{\prime}}^j\,, $$ $$ \sum
a_{ij}{X^{\prime}}^i{Y^{\prime}}^j=e'\sum a_{ij}
{Y^{\prime}}^i{X^{\prime}}^j\,, $$ that is
$a_{ij}=ea_{ij}\lambda^{i-j}$, $a_{ij}=e'a_{ji}$. Letting $i=q$
and $j=2q_0-1$ yields $e=e'=1$. Apart from the cases $i=j$;
$i=0,j=q-2q_0+1$; $i=q-2q_0+1,j=0$, some more possibilities also
arise. In fact, $i\equiv j {\pmod {q-2q_0+1}}$ together with
$i+j\leq q+2q_0-1$ does not rule out either $j=q-2_0+1+i$ for
$0<i<2q_0-1$ or $i=q-2q_0+1+j$ for $0<j<2q_0-1$. If such terms
effectively exist, then they form a polynomial of type
$G_3(X^{\prime},Y^{\prime})=\sum_{i<2q_0-1}
a_{ij}({X^{\prime}}^i{Y^{\prime}}^{q-2q_0+1+i}+{X^{\prime}}^{q-2q_0+1+i}{Y^{\prime}}^i).$
Note that $G_3(X^{\prime},Y^{\prime})$ can also be written as
$G_4(X^{\prime}Y^{\prime})({X^{\prime}}^{q-2q_0+1}+{Y^{\prime}}^{q-2q_0+1}).$
Putting $G_2(T)=G_4(T)+ T^{2q_0-1}$, we finally obtain the
required equation of $\cD$.
\end{proof}
To determine explicitly the polynomials $G$, $G_1$, and $G_2$ in
Propositions \ref{propo0} and \ref{propo1prelim} some more
computation is needed. Let $Z:=\{z \in {\mathbb F}_{q^2}\mid
b(z^q+z)=\mu^{2q_0}(\mu^q+\mu)=b^{2q_0-1}\mu^{2q_0}\}$. Let $r$ be
the line of equation $X^\prime + \lambda^2 Y^\prime=0$. As
$\xi+\lambda^2 \eta=\frac{\rho(x+\mu)}{g_0}$, by Equation
(\ref{vPx}) and Lemmas \ref{E+}, \ref{E-}, the intersection
divisor $(\pi_\ell g)^{-1}(r)$ is equal to
$$
(\pi_\ell g)^{-1}(r)=
 \left \{ \begin{array}{lll}
\sum_{z \in Z}P(z)+(2q_0+1) A_4  & \mbox{\em for
$\lambda^{q+2q_0+1}=1$} ,\\ \sum_{z \in Z,\, z \neq
\frac{\mu\lambda^q}{b},\, z\neq \frac{b+\mu\lambda^q}{b}}
P(z)+(2q_0+1) A_4 & \mbox{\em for $\lambda^{q-2q_0+1}=1$,}
                         \end{array}
          \right.
$$
where $P(z):=(1:\mu:\mu^{q_0+1}+b^{q_0}z^{q_0}:bz:b\mu
z+\mu^{2q_0+2}+b^{2q_0}z^{2q_0})$. For $P(z)$ in the support of
$(\pi_\ell g)^{-1}(r)$, the $X^\prime$-coordinate of $(\pi_\ell
g)(P(z))$ is $\lambda A(z)/B(z)$, where $
A(z):=z^{2q_0}+z+\frac{\mu^{2q_0+2}}{b^{2q_0}}$, $
B(z):=z^{q_0}+z+\frac{\mu}{b}+\frac{\mu^{q_0+1}}{b^{q_0}}$. A
straightforward computation yields $A(z)=B(z)^{2q_0}$. Hence the
affine points of $\cD\cap r$ are the points $(\lambda
B(z)^{2q_0-1},\lambda^{-1}B(z)^{2q_0-1})$, where $z$ ranges over
$Z$ for $\lambda^{q+2q_0+1}=1$, over $Z\setminus
\{\frac{\mu\lambda^q}{b},\, \frac{b+\mu\lambda^q}{b}\}$ for
$\lambda^{q-2q_0+1}=1$.
Note that the number of elements in $Z$ is $q$, but the pairwise
distinct affine points in $\cD\cap r$ are $\frac{1}{2}q$ for
$\lambda^{q+2q_0+1}=1$, $\frac{1}{2}q-1$ for
$\lambda^{q-2q_0+1}=1$. More precisely, $P(z_1)=P(z_2)$ if and
only if $B(z_1)=B(z_2)$ and this only happens when
$(z_1+z_2)^{q_0}=(z_1+z_2)$, that is either $z_1=z_2$ or
$z_1=z_2+1$  because of $g.c.d(2q_0-1,q-1)=1$. Notice also that
for $\lambda^{q-2q_0+1}=1$ and $z \in Z$, $B(z)=0$ if and only if
$z \in \{\frac{\mu\lambda^q}{b},\, \frac{b+\mu\lambda^q}{b}\}$.
Then the following lemma holds.
\begin{lemma}\label{ultimate}
With the notation of Proposition \ref{propo0}, $G$ is a polynomial
of degree $\frac{1}{2}q$ whose roots are $\{B(z)^{4q_0-2}\mid z
\in Z\}$, each of them counting once. With the notation of
Proposition \ref{propo1prelim}, $G_3$ is zero while $G_1$ is a
polynomial of degree $\frac{1}{2}q-1$ whose roots are
$\{B(z)^{4q_0-2}\mid z \in Z, B(z)\neq 0\}$, each of them counting
once.
\end{lemma}
The proposition below is the key to find the polynomials $G$ and
$G_1$.
\begin{proposition}
\label{pro1} The set $\Sigma(\lambda):=\{B(z) \ | z \in Z\}$
coincides with $Cl_0$ or $Cl_1$ according as
$\lambda^{q-2q_0+1}=1$ or $\lambda^{q+2q_0+1}=1$.
\end{proposition}
\begin{proof}
Given any $z_1\in Z$, each other element in $Z$ is written in the
form $z=z_1+c$ with $c$ ranging over $\fq$. Then
$\Sigma(\lambda)=\{B(z_1)+c^{q_0}+c \ | \ c\in \fq\}$. Hence,
either $\Sigma(\lambda)=Cl_0$ or $\Sigma(\lambda)=Cl_1$ according
as $B(z_1)\in Cl_0$ or $B(z_1)\in Cl_1$. Assume at first that
$\Sigma(\lambda)=Cl_0$. Then $z_1$ can be chosen in such a way
that $B(z_1)=0$. This implies
$z_1^{q_0}+z_1=b^{-1}\mu+b^{-q_0}\mu^{q_0+1}$. Then $z_1$ is a
solution of the system in Lemma \ref{l0}. Therefore,
$\lambda^{q^2-2q_0+1}=1$. Viceversa, for $\lambda^{q^2-2q_0+1}=1$,
the system in Lemma \ref{l0} is consistent, and taking
$(b\mu\lambda)^{-1}$ as $z_1$ $\Sigma(\lambda)=Cl_0$ follows.
Since $Cl_1=\fq \setminus Cl_0$, the proposition is proved.
\end{proof}
Let $G^+(T)\in \fq[T]$ be defined to be $$ G^+(T)=1+
\sum_{i=0}^{s-1} T^{2^iq_0}(1+T)^{2^i(q_0+1)-q_0}+T^{q/2}\,. $$
\begin{proposition}\label{roots}
With the notation of Proposition \ref{propo0}, $G(T)=eG^+(T)$,
with \mbox{$e\, \in \, {\mathbb F}_{q^4}$.}
\end{proposition}
\begin{proof}
By Proposition \ref{pro1} $B(z)$ coincides with $Cl_1$. Hence, by
Lemma \ref{ultimate}, we need to show that $G^+(a^{4q_0-2})=0 $
for any $a\in Cl_1$. Notice that $G^+(a^{4q_0-2})=0 $ if and only
if $G^+(a^{2q_0-1})=0 $, as $G^+$ is defined over ${\mathbb
F}_{2}$. Put $t=a^{2q_0-1}$. Then $t^q=t$, $t^{q_0+1}=a^{q_0}$,
$t^{2q_0+1}=a$, $t^{q/2+q_0}=a^{q/2}$. Now,
$$
(1+t)^{q_0}G^+(t)=(1+t)^{q_0}(1+t)^{q/2}+\sum_{i=0}^{s-1}(t^{q_0}(t+1)^{q_0+1})^{2i}=
$$
$$
1+t^{q_0}+t^{q/2}+t^{q/2+q_0}+\sum_{i=0}^{s-1}(t^{2q_0+1})^{2^i}+\sum_{i=0}^{s-1}(t^{q_0+1})^{2^i}+
\sum_{i=0}^{s-1}((t^{q_0})^{2^i}+(t^{2q_0})^{2^i})=
$$
$$1+t^{q_0}+t^{q/2}+t^{q/2+q_0}+\sum_{i=0}^{s-1}(t^{2q_0+1})^{2^i}+\sum_{i=0}^{s-1}(t^{q_0+1})^{2^i}+
(t^{q_0}+t^{q/2})=
$$
$$
1+a^{q/2}+\sum_{i=0}^{s-1}a^{2^i}+\sum_{i=0}^{s-1}(a^{q_0})^{2^i}=
1+a+a^2+a^4+\ldots+a^{q/2}=0.
$$
Hence the claim follows.
\end{proof}
{\bf Proof of Theorem \ref{modellong}} From Propositions
\ref{propo0} and \ref{roots} we deduce that an equation of $\cD$
is given by $eG^+(X^{\prime}Y^{\prime})=
{X^{\prime}}^{q+2q_0+1}+{Y^{\prime}}^{q+2q_0+1}$, with
$e\,\in\,{\mathbb F}_{q^4}$. Furthermore, $P:=(0,\lambda^{-1})$ is
a point of $\cD$. In fact, $P=(\pi_\ell g)
(1:\mu^{-1}:\mu^{-q_0-1}+(\lambda \mu^{-1})^{q_0}:\lambda
\mu^{-1}:\lambda \mu^{-2}+\mu^{-2q_0-2}+(\lambda
\mu^{-1})^{2q_0})$. Hence, $e=1$ and the proof is complete.

Throughout the rest of the present section we assume
$\lambda^{q-2q_0+1}=1$, and keep up the notation introduced in
Proposition \ref{propo1prelim}. Furthermore, as in Section
\ref{ordq-1}, we define $$
G^-(T)=1+\sum_{i=0}^{s-1}T^{2^i(2q_0+1)-(q_0+1)}(1+T)^{2^i}\,,\quad
\mbox{and}\quad H(T)=TG^-(T)\,. $$
\begin{lemma}\label{l7}
The roots of $H$ are the $(4q_0-2)$-{\rm th} powers of the
elements in $\Sigma(\lambda)$.
\end{lemma}
\begin{proof}
Notice that as $H$ is defined over ${\mathbb F}_{2}$ the map
$a\mapsto a^2$ is a permutation of the roots of $H$. Then the
claim follows from Proposition \ref{pro1} and Lemma \ref{l17}.
    \end{proof}
    \begin{proposition}\label{roots2}
There exists $e\in {\mathbb F}_{q^4}$ such that $G_1(T)=eG^-(T)$.
   \end{proposition}
   \begin{proof}
The claim follows from Lemmas \ref{ultimate} and \ref{l7}.
   \end{proof}
   \begin{lemma}\label{subdegree} The multiplicity of $0$ as a root of $G_2$ is $q_0-1.$
   \end{lemma}
   \begin{proof} By Lemma \ref{preimage}, $v_{P_{2}}(\xi)=q-q_0$,
$v_{P_{2}}(\eta)=1-q_0$, $v_{P_{2}}(\xi\eta)=q-2q_0+1$. Then, from
$$ v_{P_{2}}(G_1(\xi\eta))= v_{P_{2}}(G_2(\xi\eta))
+v_{P_2}(\xi^{q-2q_0+1}+\eta^{q-2q_0+1}) $$ it follows that
$v_{P_{2}}(G_2(\xi\eta))=(q-2q_0+1)(q_0-1)$, whence the claim.
     \end{proof}
\begin{proposition}\label{roots3}
There exists $e'\in {\mathbb F}_{q^4}$ such that
$G_2(T)=e'(T^{q_0-1}+ T^{2q_0-1})$.
\end{proposition}
\begin{proof}
For $\alpha\, \in \, {\bar \fq}$ root of $G_2$, let $\cC_\alpha$
be the conic of equation $X^{\prime}Y^{\prime}=\alpha$. If
$\cC_\alpha$ meets $\cD$ in a point at finite distance, say
$(x',y')$, then $G_1(x'y')=G_1(\alpha)=0$, that is $\cC_\alpha$ is
a component of $\cD$, but this is impossible.
Hence $\cC_\alpha\cap \cD$ contains no point at finite distance.
This means that $v_{P_i}(\xi\eta-\alpha)>0$ for some $i \in
\{1,2,3,4\}$. Lemma \ref{expansion} together with a
straightforward computation shows that $\alpha\in\{0,1\}$. Then
the proposition follows from Lemma \ref{subdegree}.
\end{proof}
{\bf Proof of Theorem \ref{modelsmall}} By Proposition
\ref{propo1prelim} and Lemma \ref{ultimate}, an equation of $\cD$
is
\begin{equation}
cG^-(X^{\prime}Y^{\prime})=
({X^{\prime}}^{q-2q_0+1}+{Y^{\prime}}^{q-2q_0+1})((X^{\prime}Y^{\prime})^{q_0-1}
+(X^{\prime}Y^{\prime})^{2q_0-1}),
\end{equation}
where $c\, \in \, {\mathbb F}_{q^4}$.
We prove that $c=b$. From Lemma \ref{expansion} we have for
$u=v=0$:
\[
g_0=b^{q_0-1}t+t^{q_0+1}+\ldots;\quad g_1=\mu+t+\lambda
t^{2q_0+1}+\ldots;\quad g_3=\mu+t+\lambda^{-1}t^{2q_0+1}+\ldots;
\]
whence
\[
\rho g_1/g_0=(1/b^{2q_0-1})t^{-1}(1+\ldots); \qquad \rho
g_3/g_0=(1/b^{2q_0-1})t^{-1}(1+\ldots).
\]
Since
\[
cH(\rho^2 g_1 g_3 / g_0^2)=
\rho^{q-2q_0+1}[(g_1^{q-2q_0+1}+g_2^{q-2q_0+1})/g_0^{q-2q_0+1}]
(\rho^2 g_1 g_3 / g_0^2)^{q_0} [1+(\rho^2 g_1 g_3 / g_0^2)^{q_0}]
\]
and
$$
cH(\rho^2 g_1 g_3 / g_0^2)=
c[(1/b^{2q_0-1})^2t^{-2}(1+\ldots)]^{q/2}+\ldots=
c(1/b^{2q_0-1})t^{-q}(1+\ldots),
$$
$$
\rho^{q-2q_0+1}(g_1^{q-2q_0+1}+g_2^{q-2q_0+1})/g_0^{q-2q_0+1}=
$$
\[
\rho^{q-2q_0+1}[\mu^{q-2q_0}(\lambda+\lambda^{-1})t^{2q_0+1}(1+\ldots)]/[(t^{q-2q_0+1)}
b^{(q_0-1)(q-2q_0+1)}(1+\ldots)]=
\]
\[
\rho^{q-2q_0+1}[\mu^{q-2q_0}(\lambda+\lambda^{-1})
/b^{(q_0-1)(q-2q_0+1)}] t^{-q+4q_0}(1+\ldots)=
\]
\[
b^{6-6q_0}t^{-q+4q_0}(1+\ldots),
\]
$$
(\rho^2 g_1 g_3 / g_0^2)^{q_0}=
[(1/b^{2q_0-1})^2t^{-2}(1+\ldots)]^{q_0}=
b^{2q_0-2}t^{-2q_0}(1+\ldots),
$$
$$
(1+\rho^2 g_1 g_3 / g_0^2)^{q_0}= b^{2q_0-2}t^{-2q_0}(1+\ldots),
$$
it follows that $c=b$, and the proof is complete.
\subsection{Quotient curves arising from $\D^+$}
\begin{theorem}\label{quotient+}
Let $r$ be any divisor of $q+2q_0+1$. The quotient curve
associated to the (cyclic) subgroup of order $r$ of $\D^+$ has
genus $ \frac{q_0(q-1)-1}{r}-(q_0-1)$ and is the non--singular
model of the plane $\cD_r^+$ curve of equation $$
Y^{(q+2q_0+1)/r}\Big(1+\sum_{i=0}^{s-1}
X^{2^iq_0}(1+X)^{2^i(q_0+1)-q_0}+X^{q/2}\Big)=
X^{q+2q_0+1}+Y^{2(q+2q_0+1)/r}.
$$
\end{theorem}
\begin{proof}
Let $\tau=\lambda^{(q+2q_0+1)/r}$. Then $\tau^r=1$. Let
$\varphi_r:\cD^+ \mapsto {\mathbf{P}}^2(\bar{\mathbf{F}}_q)$ be
the rational map $\varphi_r:=(1:X^{\prime}:Y^{\prime})\mapsto
(1:X^{\prime}Y^{\prime}:{Y^{\prime}}^r)$. Given a point
$Q:=(1:u:v)\in \Im(\varphi_r)$ with $v\neq 0$, let
$P:=(1:x_0:y_0)\in \varphi_r^{-1}(Q)$. For $i=1,\ldots,r$, let
$P_i:=(1:\tau^{-i}x_0:\tau^i y_0)$. Then $P_i\in \varphi^{-1}(Q)$.
On the other hand, the equation $v={Y^\prime}^r$ has exactly $r$
solutions, namely $\tau^i y_0$ with $i=1,\ldots,r$. This shows
that $\varphi_r^{-1}(Q)=\{(1:\tau^{-i}u/y_0:\tau
y_0)|i=1,\ldots,r\}$. In particular, $\varphi_{r}$ has degree $r$.
By Corollary \ref{c1}, it turns out that the non--singular model
of $\varphi_r(\cD^+)$ is the quotient curve of $\cX$ with respect
to the automorphism $[B]^{(q+2q_0+1)/r}$ of $\cX$. A
straightforward computation gives the desired equation.
\end{proof}
\subsection{Quotient curves arising from $\D^-$}
\begin{theorem}\label{quotient-}
Let $r$ be any divisor of $q-2q_0+1$. The quotient curve
associated to the (cyclic) subgroup of order $r$ of $\D^-$ has
genus $ \frac{q_0(q-1)+1}{r}-(q_0+1)$ and is the non--singular
model of the plane $\cD_r^-$ curve of equation
\begin{equation}
bY^{\frac{q-2q_0+1}{r}}\Big(1+\sum_{i=0}^{s-1}X^{2^i(2q_0+1)-q_0-1)}(1+X)^{2^i}\Big)
=(X^{q-2q_0+1}+Y^{\frac{2(q-2q_0+1)}{r}})(X^{q_0-1}+X^{2q_0-1}).
\end{equation}
\end{theorem}
\begin{proof}
The proof is similar to the proof of Theorem \ref{quotient+}.
\end{proof}

\section{Quotient curves arising from non-tame subgroups}\label{nontame}

In the following sections we will investigate the quotient curves
of the DLS-curve arising from its automorphism groups of even
order. We will give a method for computing  the genera of such
curves, and in several cases we will also provide an equation for
them. Our approach is similar to that employed in \cite{g-sti-x}
and \cite{aq} where the function field point of view was used to
investigate the analogous problem for the Hermitian curve. In that
context, $\fq$-automorphisms are viewed as elements of the
automorphism group of the function field. In our case,
$\fq(\cX)=\fq(x,y)$ with $x^{2q_0}(x^q+x)=by^q+by$, and
$\Aut(\cX)\cong {\cS}z(q)$, see Remark \ref{remark2}. The
extension  $\fq(\cX)|\fq(x)$ is Galois of degree $q$, and $x$ has
a unique pole in $\fq(\cX)$ that we denote by $\cP_{\infty}$. Such
a place is totally ramified in $\fq(\cX)$, while all other
rational places of $\fq(x)$ split completely in $\fq(\cX)|\fq(x)$.
The Galois group of $\fq(\cX)|\fq(x)$ is $\bar{{\bf
T}}_0:=\{\psi_{0,c}\mid c\in\fq\}$ with $\psi_{a,c}$ as in Section
\ref{dls}. Note that $\bar{{\bf T}}_0$ comprises the identity and
the elements of order $2$ of the Sylow $2$-subgroup ${\bar{\bf
T}}=\{\psi_{a,c}\mid a,c\in\fq \}$ of $\Aut(\cX)$. Let $\cX_{\cU}$
be the quotient curve of $\cX$ associated to a subgroup $\cU$ of
$\Aut(\cX)$ of even order. In computing the genus $g_{\cU}$ of
$\cX_{\cU}$ by means of the Hurwitz genus formula, the essential
problem is to compute $\deg(\Diff(\fq(\cX)|\fq(\cX_{\cU}))$. Since
$\fq(\cX)|\fq(\cX_{\cU})$ is a non-tame extension, knowing the
order of $\cU$ and its action on places of $\fq(\cX)$ is not
sufficient to compute $\deg(\Diff(\fq(\cX)|\fq(\cX_{\cU}))$.
However as in \cite{g-sti-x}, the Hilbert different's formula see
\cite[Prop.III.5.12, Theor. III.8.8]{sti} allows us to overcome
this difficulty. Let $\mathbb P$ denote the set of all places of
$\fq(\cX)$. For $\cP\, \in \, {\mathbb P}$, the Hilbert
different's formula states that the different exponent $d(\cP)$ of
$\cP$ with respect to the extension $\fq(\cX)|\fq(\cX_{\cU})$ is
$$
d(\cP)=\sum_{\sigma\in \cU\setminus \{1\}, \sigma(\cP)=\cP}
i_\cP(\sigma),
$$
where $i_\cP(\sigma)= v_P(\sigma(t)-t)$, $t$ being a local
parameter at $\cP$. Hence,
$$
\deg (\Diff (\fq(\cX)|\fq(\cX_\cU))=\sum_{1\neq\sigma\in \cU}\Big(
\sum_{ \cP\in {\mathbb P} , \sigma(\cP)=\cP}
 i_\cP(\sigma)\Big) \,.
$$
\begin{proposition}
For $1\neq \psi_{a,c}\, \in \,{\bar {\bf T}}$,
\begin{equation}\label{order2}
i_{\cP_\infty}(\psi_{a,c})=
\begin{cases}
2q_0+2 & \text{for $a=0$}\, ,\\ 2 & \text{for a $\neq 0$}\, .
\end{cases}
\end{equation}
\end{proposition}
\begin{proof}
Let $t:=f_3/f_4$. Then $t$ is a local parameter at $\cP_\infty$.
By straightforward computation $$
\psi_{a,c}(t)=\frac{a^{2q_0}x+by+c}{f_4+\alpha}, $$ where
$\alpha=a^{2q_0+1}x+cx+aby+a^{2q_0+2}+ac+c^{2q_0}$. Hence $$
\psi_{a,d}(t)-t=\frac{xf_4a^{2q_0}+cf_4+by\alpha}{f_4(f_4+\alpha)}\,.
$$ Then the assertion follows from the following computation. Let
$f_1:=x,\, f_2:=x^{q_0+1}+b^{q_0}y^{q_0}, \, f_3:=by,\,
f_4:=bxy+x^{2q_0+2}+b^{2q_0}y^{2q_0}$. Then
\begin{enumerate}
\label{lem1}
\item $v_{\cP_\infty}(f_1)=-q$,
\item $v_{\cP_\infty}(f_3)=-2q_0-q$,
\item $v_{\cP_\infty}(f_2)=-q_0-q$,
\item $v_{\cP_\infty}(f_4)=-2q_0-q-1$.
\end{enumerate}
In fact, we have
\begin{enumerate}
\item $\cP_\infty$ is
the only pole of  $x$ and  $[\fq (\cX):\fq (x)]=q$.
\item $\cP_\infty$ is
the only pole of $by$  and $[\fq (\cX):\fq (by)]=2q_0+q$ .
\item $f_2^{2q_0}=by^q+x^{q+2q_0}=by+x^{2q_0+1}$ and
$v_{\cP_\infty}(by+x^{2q_0+1})=-q(2q_0+1)$.
\item $f_4^{q_0}=b^{q_0}x^{q_0}y^{q_0}+x^{q+2q_0}+by^q=
b^{q_0}x^{q_0}(y^{q_0}+x^{q_0 +1})+by$ whence
$q_0v_{\cP_\infty}(f_4)=-(qq_0+q+q_0)$.
\end{enumerate}
\end{proof}
\begin{corollary}[Hilbert different's formula for the
DLS-curve] \label{HilDLS} For a subgroup $\cU$ of $\Aut(\cX)$ let
$N_i$ denote the number of elements of $\cU$ of order $i$ with
$i>1$. Then $$\deg (\Diff (\fq(\cX)|\fq(\cX_\cU))=\sum_{i>1}
k_iN_i$$ where
\begin{equation}
\label{formula}
 k_i \ = \left \{ \begin{array}{lll}
         2q_0+2       & \mbox{\em for \ $i=2$ \/}, \\
         2            & \mbox{\em for \ $i=4$ \ and  \ for \ $i\mid q-1$ \/}, \\
         4            & \mbox{\em for \ $i\mid q-2q_0+1$ \/}, \\
         0            & \mbox{\em otherwise. \/}
                         \end{array}
          \right.
\end{equation}
\end{corollary}
\begin{proof} A full set of conjugacy class representatives
of non trivial elements in ${\cS}z(q)$ is $\{[T_{0,1}],[T_{1,1}],
[N_d], ({{\bf g}^+})^j, ({{\bf g}^-)}^j\}$ where $d$ ranges over
$\fq\setminus \{0\}$, ${\bf g}^\pm$ generates ${\bf D}^\pm$ and
$j=1,\ldots q\pm2q_0$. For the corresponding automorphisms in
$\Aut(\cX)$, (\ref{formula}) for $i=2,4$ comes from
(\ref{order2}), while (\ref{formula}) for odd $i$ follows from the
fact that the number of fixed points on $\cX$ of $[N_d]$, $({{\bf
g}^+})^j$, and $({{\bf g}^-})^j$ is equal to $2,0$ and $4$,
respectively. Finally, (\ref{formula}) holds true for any
non-trivial element of $\Aut(\cX)$ since elements in ${\cS}z(q)$
(and hence in $\Aut(\cX)$) of the same order are pairwise
conjugate.
\end{proof}

    \section{Quotient curves of $\cX$ associated to $2$ subgroups}
\label{2subgroups}

Throughout this section we use the following notation:
\begin{itemize}
\item $\cU$ is a subgroup of ${\bar {\bf T}}$;
\item $\cU_2$ is the subgroup of $\cU$ consisting of all elements of order $2$ together with the
identity;
\item $\cX_\cU$ is the quotient curve of $\cX$ arising from $\cU$;
\item $g_{\cU}$ is the genus of $\cX_\cU$.
\end{itemize}
We begin by giving a formula for computing the genus.
\begin{proposition}\label{2sbgrp}
Let $\cU$ have order $2^u$. If $\cU_2$ has order $2^v$, then
$$g_\cU=2^{s-u+v}(2^{2s+1-v}-1).$$
\end{proposition}
\begin{proof}
By Corollary \ref{HilDLS}
$$
d(\cP_\infty)=2(2^{u}-2^v)+(2q_0+2)(2^v-1)\,.
$$
Then the Hurwitz genus formula yields $$
2q_0(q-1)-2=2^{u}(2g_\cU-2)+2(2^{u}-2^v)+(2q_0+2)(2^v-1)\,, $$
that is $$ g_\cU=q_0(q-2^v)/2^{u}=2^{s-u+v}(2^{2s+1-v}-1)\,. $$
\end{proof}
Proposition \ref{2sbgrp} rises the problem of classifying the
subgroups of $\bar {\bf T}$ in terms of the number of their
elements of order $2$. Such a general problem is computationally
beyond our possibility, because $\bar{\bf T}$ contains a huge
number of pairwise non-conjugate subgroups. What we do here is to
prove some results which are useful to investigate special cases.
The following proposition states some numerical conditions on $u$
and $v$.
\begin{proposition}
\label{legr1}
Let $\cU$ have order $2^u$. If $\cU_2$ has order $2^v$ then
\begin{itemize}
\item[I)]   $u\leq 2v$;
\item[II)]  for every integer $u^{\prime}$ with $v \leq u^{\prime}\leq u$
there is a subgroup of $\cU$ of order $2^{u^{\prime}}$;
\item[III)] for $2^v<q$ we have $u-v\le s$.
\end{itemize}
\end{proposition}
\begin{proof}
The map $\Phi:{\bar{\bf T}}\to\fq$ given by $\Phi({\psi}_{a,c})=a$
is a homomorphism from $\bar{\bf T}$ onto the additive subgroup of
$\fq$. The restriction of $\Phi$ to $\cU$ is the homomorphism
$\Phi_{|\cU}$ with kernel $\Ker(\Phi_{|\cU})=\{\psi_{0,c}|
c\in\fq\}$ isomorphic to $\cU_2$. Both $\Im(\Phi_{|\cU})$ and
$\Ker(\Phi_{|\cU})$ are linear subspaces of $\fq$ regarded as a
vector space over ${\mathbb{F}}_2$. This yields $u=v+w$ where
$2^w$ denotes the order of $\Im(\Phi_{|\cU})$ viewed as an
additive subgroup of $\fq$. Now, since
$\psi_{a,c}^2=\psi_{0,a^{2q_0+1}}$ holds, we have that
$\Im(\Phi_{|\cU})^{2q_0+1}:=\{a^{2q_0+1}|a\in \Im(\Phi_{|\cU})\}$
is a subset of $\Ker(\Phi_{|\cU})$. As $a\mapsto a^{2q_0+1}$ is
one-to-one map of $\fq$, we have $w\leq v$, whence assertion I)
follows. We note that the factor--group $\cU/\cU_2$ is an
elementary abelian of order $2^w$ because $\cU_2$ contains all
elements of $\cU$ of order $2$. Assertion II) follows from the
well known fact that the converse of the Lagrange theorem holds
for any elementary abelian group. Proposition \ref{2sbgrp}
together with the fact that $g_\cU$ must be an integer gives
assertion III).
\end{proof}
In the case where $\cU=\bar{\bf T}$ and $\cU_2=Z(\bar{\bf T})$,
assertion II) in Proposition \ref{legr1} has the following
corollary.
\begin{lemma}
\label{legr2} For every integer $u$ with $s\leq u \leq 2s+1$, there is a
subgroup $\cU$ of $\bar{\bf T}$ of order $2^u$ containing all
elements of $\bar{\bf T}$ of order $2$.
\end{lemma}
The existence of a subgroup of $\bar{\bf T}$ with a given number
of elements of order $2$ is ensured by the following lemmas.
\begin{lemma}
\label{legr3} Let $\cH$ be an elementary abelian subgroup of
$\bar{\bf T}$ of order $2^v$. Then there is a subgroup $\cU$ of order  $2^{v+1}$
such that $\cU_2$ coincides $\cH$.
\end{lemma}
\begin{proof}
Since the normaliser of the Sylow $2$-subgroup $\bf T$ of
${\cS}z(q)$ acts transitively on the set of elements of order $2$
in $\bf T$, we may assume $\psi_{0,1}\in \cH$. Then the group
generated by $\cH$ and $\psi_{1,0}$ has the required property.
\end{proof}

\begin{lemma}\label{legr4}
Let $\cB$ be an additive subgroup of $\fq$ of order $2^v$ with
$0\le v \le 2s+1$. If there exists an additive subgroups $\cA$ of
$\fq$ of order $2^{u-v}$ such that $\cA^{2q_0+1}\subseteq \cB$,
then there is a subgroup $\cU$ of $\bar{\bf T}$ of order $2^u$
such that $\Ker(\Phi_{|\cU})=\cB$ and $\Im(\Phi_{|\cU})= \cA$.
\end{lemma}
\begin{proof} We keep the notation introduced in the proof of
Proposition \ref{legr1}. The proof of Lemma \ref{legr4} is by
induction on $u$. First we consider the case $u=v$. We have
$\cA=\{0\}$, and hence $\cA^{2q_0+1}=\{0\}$. Let
$\cU=\{\psi_{0,c}|c\in \cB\}$. Then $\cU$ has the required
properties. Suppose now that $u>v$. As $\cA$ is an elementary
abelian group, it contains a subgroup $\cA_0$ of index $2$, that
is of order $2^{u-1-v}$. Since $\cA_0^{2q_0+1}\subseteq
\cA^{2q_0+1}$, there is by induction a subgroup $\cU_0$ in
$\bar{\bf T}$ with $\Ker(\Phi_{|\cU_0})=\cB$ and
$\Im(\Phi_{|\cU_0})= \cA_0$. For a fixed $\beta=\psi_{a,c}\in
\bar{\bf T}$ with $a\in \cA\setminus \cA_0$, let $\cU$ be the
subgroup of $\bar{\bf T}$ generated by $\cU_0$ together with
$\beta$. We show that $\cU$ has order $2^{u-v}$, that is
$\cU=\beta U_0 \cup U_0$. To do this, it is enough to check that
$\beta U_0=U_0 \beta$. For every element $\gamma:=\psi_{a_0,c_0}$
of $\cU_0$, we have $(\beta\gamma)^2=\psi_{0,(a+a_0)^{2q_0+1}}$.
Since $a+a_0\in \cA$, we have $(a+a_0)^{2q_0+1}\in \cB$. As
$\Ker(\Phi_{|\cU_0})=\cB$, we obtain indeed
$(\beta\gamma)^2\in\cU_0$. Using this fact together with two more
properties, namely that $\bf T$ has exponent $4$ and that every
involutory element in $\bf T$ is in the center $Z(\bf T)$, we have
$$ \beta\gamma=\beta \gamma (\beta
\gamma^4\beta^3)=(\beta\gamma)^2\gamma^3\beta^3=\gamma^3(\beta\gamma)^2\beta^2\beta\in
\cU_0\cU_0\cU_0\beta=\cU_0\beta. $$ Finally,
$\Phi(\beta\gamma)=\Phi(\beta)+\Phi(\gamma)=a+a_0$ implies not
only that $\Im(\Phi_{|\cU})=\cA$ but also that no element in
$\beta \cU_0$ is in $\Ker(\Phi_{|\cU})$, whence
$\Ker(\Phi_{|\cU})=\Ker(\Phi_{|\cU_0})=\cB$ follows.
\end{proof}
\begin{lemma}\label{legr5}
For two positive integers $u,v$ with $u\geq v$, there is a
subgroup $\cU$ of order $2^u$ such that $\cU_2$ has order $2^v$,
provided that one of the following holds:
\begin{itemize}
\item  $v\le 2s+1$ and $u \leq v+\log_2{(v+1)}$,
\item  $(u-v)|(2s+1)$ and $v\ge u-v$.
\end{itemize}
\end{lemma}
\begin{proof}
For any additive subgroup $\cA$ of $\fq$ of order $2^{u-v}$, the
additive subgroup $\cB^{\prime}$ of $\fq$ generated by all
elements in $\cA^{2q_0+1}$ has order at most $2^{2^{u-v}-1}$. In
fact, $\fq$ can be viewed as a vector space over its subfield
$\mathbb{F}_2$, and the subspace generated by $\cA^{2q_0+1}$ has
dimension at most $2^{u-v}-1$. Suppose at first that both $v\le
2s+1$ and $u\le v+\log_2{(v+1)}$ hold. Then there exists an
additive subgroup $\cB$ of $\fq$ of order $2^v$ containing
$\cB^{\prime}$, and the first claim follows from Lemma 7.5. Now
suppose that $v\ge u-v\ge 0$, $(u-v)|(2s+1)$. Then there exists a
subfield of $\mathbb{F}_m$ of $\fq$ of order $m=2^{u-v}$. Let
$\cB$ be any additive subgroup of order $2^v$ containing the
additive group $\cA$ of $\mathbb{F}_m$. Again Lemma 7.5 proves the
claim.
\end{proof}

\begin{remark}\label{rem2}
Lemma \ref{legr3} does not hold true for subgroups $\cU$ of order
$2^{v+\ell}$ with $\ell>1$, as the following example shows. Fix an
element $e\in{\mathbb{F}}_q \backslash {\mathbb{F}}_2$. The set
$\cU_2=\{\psi_{0,0},\psi_{0,1},\psi_{0,e},\psi_{0,e+1}\}$ is an
elementary abelian group of order $2^v$ with $v=2$. Assume that
there is a subgroup $\cU$ of $\bar{\bf T}$ of order $2^4$ whose
elements of order $2$ are those of $\cU_2$. Then there are three
pairwise distinct non-zero elements
$a_1=1,a_2,a_3\in{\mathbb{F}}_q$ and three elements
$c_1,c_2,c_3\in{\mathbb{F}}_q$ such that $\psi_{a_i,c_i}$,
$i=1,2,3$ together with $\psi_{0,0}$ form a complete set of
representatives of the cosets of $\cU/\cU_2$. Furthermore, $\Im
(\Phi)=\{0,1,a_2,a_3\}$. Hence $a_3=1+a_2$. On the other hand,
$\{a_2^{2q_0+1},a_3^{2q_0+1}\}=\{e,e+1\}$ by the proof of Lemma
\ref{legr3}. Thus $a_3^{2q_0+1}=1+a_2^{2q_0+1}$. From these
results, $(1+a_2)^{2q_0+1}=1+a_2^{2q_0+1}$. Hence
$a_2^{2q_0-1}=1$. But this is impossible as $a_2\neq 1$ and
$2q_0-1$ is coprime to $q-1$.
\end{remark}

In some cases we are able to provide an equation for the quotient
curve $\cX_\cU$.
\begin{theorem}\label{subfield}
For a subfield ${\mathbb F}_{q^{\prime}}$ of $\fq$, let $\cU$ be
the elementary abelian subgroup of $\bar{\bf T}$ consisting of all
automorphisms $\psi_{0,c}$ with $b^{-1}c\in{\mathbb
F}_{q^{\prime}}$. Then the quotient curve $\cX_\cU$ has genus
$g_{\cU}=q_0(q/q^{\prime}-1)$ and is a non-singular model over
$\fq$ of the irreducible plane curve of equation
\begin{equation}
X^{2q_0}(X^q+X)=b\sum_{i=0}^{n-1}Y^{(q')^i}
\end{equation}
where $q=(q')^n$.
\end{theorem}
\begin{proof}
Let  $\Phi:\cC_b\rightarrow \P^2(\bfq)$ be the rational map
$\Phi:(1:X:Y)\mapsto (1:X:Y^{q'}+Y)$. Given a point $Q:=(1:u:v)\,
\in \, \Im (\Phi)$, let $P:=(1:x_0:y_0)\, \in\, \Phi^{-1}(Q)$. For
$\omega \in \{b^{-1}c\in \fq \mid \psi_{0,c}\in \cU\}={\mathbb
F}_{q'}$, let $P_i:=(1:x_0:y_0+\omega)$. Then $P_i\, \in \,
\Phi^{-1}(Q)$. On the other hand, the equation $\eta=Y^{q'}+Y$ has
exactly $q'$ solutions, namely $\{y_0+\omega \mid \omega \in
{\mathbb F}_{q'}\}$. This shows that
$\Phi^{-1}(Q)=\{(1:x_0:y_0+b^{-1}c) \mid \psi_{0,c} \in \cU\}$.
Hence, the non-singular model of $\Phi(\cC_b)$ is the quotient
curve of $\cX$ with respect to $\cU$. Then a straightforward
computation showing that $$
y^q+y=(y^{q'}+y)+(y^{q'}+y)^{q'}+\ldots (y^{q'}+y)^{(q')^{n-1}} $$
completes the proof.
\end{proof}
\begin{theorem}
\label{ord4} For a cyclic subgroup $\cU$ of order $4$ of
$\Aut(\cX)$, the quotient curve $\cX_\cU$ of $\cX$ associated to
$\cU$ has genus $g_{\cU}=\frac{1}{4}q_0(q-2)$ and it is a
non-singular model over $\fq$ of the irreducible plane curve of
equation
\begin{equation}
\sum_{i=0}^{2s}X^{2^i}+\sum_{i=0}^{s}X^{2^i}\Big( \sum_{j=i}^s
X^{2^j}\Big) +\sum_{i=s+1}^{2s}
X^{2^i}\Big(\sum_{j=0}^{i-s-2}X^{2^j}\Big)^{2q_0}= \sum
_{i=0}^{2s}Y^{2^i}\,.
\end{equation}
\end{theorem}
\begin{proof}
Since the cyclic subgroups of ${\cS}z(q)$ of order $4$ are
pairwise conjugate under ${\cS}z(q)$, we may assume $\cU$ to be
generated by $\psi_{1,0}$. Let  $\Phi:\cC_b\rightarrow \P^2(\bfq)$
be the rational map $\Phi:(1:X:Y)\mapsto
(1:X^2+X:b^2Y^2+bY+X^3+X)$. Given a point $Q:=(1:u:v)\, \in \, \Im
(\Phi)$, let $P:=(1:x_0:y_0)\, \in\, \Phi^{-1}(Q)$. Then it is
easily seen that $\Phi^{-1}(Q)=\{\psi_{1,0}^i(P) \mid
i=1,\ldots,4\}$. Hence a non-singular model of $\Phi(\cC_b)$ is
the quotient curve of $\cX$ arising from $\cU$. Now, since $$
\sum_{i=0}^{2s}(b^2Y^2+bY+X+X^3)^{2^i}=b(Y^q+Y)+\sum_{i=0}^{2s}(X+X^3)^{2^i}\,,
$$ we only have to show that $$
\sum_{i=0}^{2s}(X+X^3)^{2^i}+X^{2q_0}(X^q+X)= $$ $$
\sum_{i=0}^{2s}(X+X^2)^{2^i}+\sum_{i=0}^{s}(X+X^2)^{2^i}\Big(
\sum_{j=i}^s (X+X^2)^{2^j}\Big) +\sum_{i=1}^{s-1}\Big(
(X+X^2)^{2^i} (\sum_{j=0}^{i-1}(X+X^2)^{2^j})\Big)^{2q_0}. $$ This
follows from the following two equations: $$
\sum_{i=0}^{2s}(X+X^3)^{2^i}+X^{2q_0}(X^q+X)=
\sum_{i=0}^{2s}(X+X^2)^{2^i}+ \sum_{i=0}^{2s}(X(X+X^2))^{2^i}+
X^{2q_0}\sum_{i=0}^{2s}(X+X^2)^{2^i}$$ $$
=\sum_{i=0}^{2s}(X+X^2)^{2^i}+\sum_{i=0}^{2s}(X^{2^i}+X^{2q_0})(X+X^2)^{2^i}\,,
$$ $$ X^{2^i}+X^{2q_0}=\begin{cases} \sum_{j=i}^{s}(X+X^2)^{2^j} &
\text{if \quad $i<s+1$ }\, ,\\ 0 & \text{if \quad $i=s+1$ }\, ,\\
(\sum_{j=0}^{i-s-2}(X+X^2)^{2^j})^{2q_0} & \text{if \quad $i>s+1$
}\, .
\end{cases}
$$ \end{proof}

    \section{Quotient curves arising from subgroups of order $2^ur$
fixing a place; $u>1$ and $r>1$ is a divisor of $q-1$}\label{rdueu}

We keep the notation introduced in Section \ref{2subgroups}. In
addition, let $\bar{\bf N}=\{\gamma_{d}\mid d\in \fq\}$ with
$\gamma_{d}$ as in (\ref{gamma}). Then ${\bar{\bf T}}{\bar {\bf
N}}$ is the normaliser of $\bar{\bf T}$. Assume that ${\bar{\bf
T}}{\bar {\bf N}}$ contains a subgroup $\cU$ of order $2^ur$ with
$r\mid q-1$. Every subgroup of order $2^ur$ is conjugate to $\cU$
under $\Aut(\cX)$. Note that ${\bar{\bf T}}{\bar {\bf N}}$ viewed
as a permutation group on the set of all $\fq$-rational places
different from $\cP_{\infty}$, is a Frobenius group with kernel
$\bar{\bf T}$ and nucleus $\bar{\bf N}$. Thus, the order of a
non-trivial element $\sigma \in {\bar{\bf T}}{\bar {\bf N}}$ is
either a $2$-power or a divisor of $q-1$ according as $\sigma$
belongs to ${\bar{\bf T}}$ or does not. This together with
Corollary \ref{HilDLS} gives the following result.
\begin{theorem}
For a subgroup $\cU$ of $\bar{\bf T}\bar{\bf N}$ of order $2^ur$,
let the subgroup $\cU_2$ of $\cU$ consist of all elements of order
$2$ together with the identity. If $\cU_2$ has order $2^v$, then
$$ g_{\cU}=\frac{1}{r}[2^{s-u+v}(2^{2s+1-v}-1)]\,. $$
\end{theorem}
\begin{proof}
By the Hurwitz genus formula, $$
2q_0(q-1)-2=r2^u(2g_\cU-2)+2(2^u-2^v)+(2q_0+2)(2^v-1)+2(
r2^u-2^u)$$
\end{proof}

\section{Quotient curves arising from dihedral subgroups of order $2r$ with
 a divisor $r>1$ of $q-1$}\label{rdivq-1}

The normaliser $N_{\Aut(\cX)}(\bar{\bf {N}})$ of $\bar{\bf N}$ in
is the dihedral group of order $2(q-1)$ which comprises $\bar{\bf
N}$ together with a coset consisting entirely of elements of order
$2$. Let $\cU$ be a subgroup of $\Aut(\cX)$ of order $2r$ with a
divisor $r>1$ of $q-1$. Up to conjugacy, $\cU$ is a subgroup of
$N_{\Aut(\cX)}({\bar{\bf N}})$. Hence $\cU$ has $r-1$ non-trivial
elements of odd order and each of the remaining $r$ elements in
$\cU$ have order $2$. The argument in Section \ref{rdueu}
depending on the Hurwitz genus and the Hilbert different formulas,
enable us to compute the genus of the quotient curve $\cX_{\cU}$
of $\cX$ arising from $\cU$. To find an equation for $\cX_{\cU}$
we also need the Waring formula in two indeterminates, say $U$ and
$V$, see \cite[Theorem 1.76]{LN}:
\begin{result}
\label{Ris6}
\begin{equation}
\label{waring} U^k+V^k=\sum
(-1)^{i+j}\frac{(i+j-1)!k}{i!j!}(U+V)^i(UV)^j
\end{equation}
where the summation is extended over all pairs $(i,j)$ of
non--negative integers for which $i+2j=k$ holds.
\end{result}
\begin{theorem}
\label{teordivq-1} Let $\cU$ be a subgroup of $\Aut(\cX)$ of order
$2r$ with $r\mid q-1$, and $r>1$. Then the quotient curve
$\cX_\cU$ has genus $g_\cU=q_0(q-r-1)/2r$ and it is a
non--singular model over $\fq$ of the irreducible plane curve of
equation is
$$1+\sum_{i=0}^{s-1}X^{2^i(2q_0+1)-(q_0+1)}(1+X)^{2^i}=\sum
(-1)^{i+j}\frac{(i+j-1)!k}{i!j!}Y^i(X^{rj}(X^{q_0}+1),$$ where the
summation is extended over all pairs $(i,j)$ of non--negative
integers with $i+2j=(q+2q_0+1)/r$.
\end{theorem}
\begin{proof}
With the notation introduced in Section \ref{ordq-1}, let
$\phi_r:\pi_{\ell}(\cX)\mapsto \P^2(\bfq)$ be the rational map
defined as $\phi_r:(1:X^{\prime}:Y^{\prime})\mapsto
(1:X^{\prime}Y^{\prime}:{X^{\prime}}^r+{Y^{\prime}}^r)$. We argue
as in the proof of Theorem \ref{ff2}. Given a point $Q:=(1:u:v)\in
\Im (\phi_r)$ with $u,v\neq 0$, let
$P:=(1:x_0:y_0)\in\phi_r^{-1}(Q)$, For $i=1,\ldots,r$, let
$P_i:=(1:\tau^{-q_0i}x_0:\tau^{q_0i}y_0),
P_i^{\prime}:=(1:\tau^{q_0i}y_0:\tau^{-q_0i}x_0)$ with an element
$\tau$ of order $(q-1)/r$ in the multiplicative group of $\fq$.
Then both $P_i$ and $P_i^{\prime}$ are in $\phi^{-1}(Q)$. On the
other hand, the equation ${Y^{\prime}}^{2r}+v{Y^{\prime}}^r+u^r=0$
has $2r$ pairwise distinct solutions. Hence, $\phi_r$ has degree
$2r$. So, the quotient curve $\cU_{\cX}$ of $\cX$ associated to
$\cU$ is the non singular model over $\fq$ of the irreducible
plane curve $\phi_r\pi_{\ell}(\cX)$. The equation of the latter
curve derives from (\ref{gminus}) taking into account Result
\ref{Ris6} applied to $U=X^r$, $V=Y^r$ and $k=(q-1)/r$.
\end{proof}

\section{Quotient curves arising from subgroups of order $2r$ with
a divisor $r>1$ of $q\pm 2q_0+1$}

The long as well as the short Singer subgroup of ${\cS}z(q)$ is the maximal cyclic subgroup of a dihedral group of ${\cS}z(q)$. Up to conjugacy under ${\cS}z(q)$,
such a dihedral group $\cU$ comprises the Singer subgroup $\bf D$
generated by $[B]^{(q\pm 2q_0+1)/r}$ together with the coset
$[W]\bf D$ of elements of order $2$. The statements in Section
\ref{rdivq-1} hold true if $q-1$ is replaced by $q+2q_0+1$ for the
long Singer subgroup and by $q-2q_0+1$ for the short Singer
subgroup.
\begin{theorem}\label{dihedral+}
Let $r>1$ be a divisor of $q+2q_0+1$. The quotient curve
$\cX_{\cU}$ of $\cX$ associated to a subgroup $\cU$ of order $2r$
has genus
$$g_{\cU}=\frac{1}{2}\left[\frac{q_0(q-1)-1}{r}-(q_0-1)\right].$$
Furthermore, $\cX_{\cU}$ is the non--singular model over
${\mathbb{F}}_{q^4}$ of the irreducible plane curve of equation
\begin{equation}
1+\sum_{i=0}^{s-1} X^{2^iq_0}(1+X)^{2^i(q_0+1)-q_0}+X^{q/2}=\sum
(-1)^{i+j}\frac{(i+j-1)!}{i!j!}X^{ri}Y^j,
\end{equation}
where the summation is extended over all pairs $(i,j)$ of
non--negative integers with $i+2j=(q+2q_0+1)/r$.
\end{theorem}
\begin{proof}
To find the equation we will argue as in the proof of Theorem
\ref{teordivq-1}. Let $\psi_r:\cD^+ \mapsto
{\mathbf{P}}^2(\bar{\mathbf{F}}_q)$ be the rational map
$\psi_r:=(1:X:Y)\mapsto (1:XY:X^r+Y^r)$. Given a point
$Q:=(1:u:v)\in \Im(\psi_r)$ with $u\neq 0,v\neq 0$, let
$P:=(1:x_0:y_0)\in \psi_r^{-1}(Q)$. For $i=1,\ldots,r$, let
$P_i:=(1:\tau^{-i}x_0:\tau^i y_0),R_i:=(1:\tau^i
y_0:\tau^{-i}x_0)$, with $\tau$ an element of order $(q-1)/r$ in
the multiplicative group of $\fq$. Then both $P_i$ and $R_i$  are
in $\psi_r(Q)^{-1}$. On the other hand, if
$P:=(1,u',v')\in\psi_r(Q)^{-1}$, then $v'=u (u')^{-1}$, and hence
$(u')^{2r}+v (u')^r+u=0$. Since the latter equation has $2r$
pairwise distinct solutions in $u'\in \bar{\mathbb{F}}_q$, we
obtain that $\psi_r(Q)^{-1}=\{P_i|i=1,\ldots,r\}\cup$
$\{R_i|i=1,\ldots,r\}.$ This shows that $\psi_r$ has degree $2r$.
By Corollaries \ref{c1} and \ref{c2}, the non--singular model of
$\psi_r(\cD^+)$ is the quotient curve $\cX_{\cU}$ of $\cX$ arising
from $\cU$. The computation for the equation of $\psi_r(\cD^+)$
can be carried out as before, by applying (\ref{waring}) for
$U=X^r,V=Y^r,$ and $k=(q+2q_0+1)/r.$
\end{proof}
\begin{theorem}\label{dihedral-} Let $r>1$ be a divisor of
$q-2q_0+1$. The quotient curve $\cX_{\cU}$ associated to a
subgroup $\cU$ of order $2r$ has genus $$g_{\cU}=\frac{1}{2}\left[
\frac{q_0(q-1)+1}{r}-(q_0+1)\right].$$ Furthermore, $\cX_{\cU}$ is
the non--singular model over ${\mathbb{F}}_{q^4}$ of the
irreducible plane curve of equation
$$b\Big(1+\sum_{i=0}^{s-1}X^{2^i(2q_0+1)-(q_0+1)}(1+X)^{2^i}\Big)
=(X^{q_0-1}+X^{2q_0-1})\sum
(-1)^{i+j}\frac{(i+j-1)!}{i!j!}X^{ri}Y^j,$$
where the summation is extended over all pairs $(i,j)$ of
non--negative integers with $i+2j=(q+2q_0+1)/r$.
\end{theorem}
\begin{proof} Let $\psi_r:\cD^+ \mapsto
{\mathbf{P}}^2(\bar{\mathbf{F}}_q)$ be the rational map
$\psi_r:=(1:X:Y)\mapsto (1:XY:X^r+Y^r)$. Arguing as in the proof
of Theorem \ref{dihedral-}, it turns out that a non-singular model
of $\psi_r(\cD^+)$ is the quotient curve of $\cX$ arising from
$\cU$. Again, the computation for the equation of $\psi_r(\cD^-)$
can be carried out as before, by applying (\ref{waring}) for
$U=X^r,V=Y^r,$ and $k=(q-2q_0+1)/r.$
\end{proof}

\section{Quotient curves arising from subgroups of order $4r$ with
a divisor $r>1$ of $q\pm 2q_0+1$}

The normaliser $N_{{\cS}z(q)}(\bf D^+)$ of the long Singer
subgroup $\bf D^+$ of ${\cS}z(q)$ is a Frobenius group of order
$4(q+2q_0+1)$ with kernel $\bf D^+$ and complement $C_4$ where
$C_4$ is the cyclic group generated by the linear collineation
associated to the matrix
$$ \left( \begin{array}{ccccc} 1 & 0 & 0 & 1 & 1 \\ 1 & 1 & 0 & 1
& 0 \\ 0 & 0 & 1 & 0 & 0 \\ 0 & 0 & 1 & 0 & 0 \\ 1 & 1 & 0 & 0 & 1
\end{array} \right)\,.$$
Let $\cU$ be a subgroup of order $4r$ such that $r>1$ divides
$q+2q_0+1$. Up to conjugacy under ${\cS}z(q)$, $\cU$ is a subgroup
of $N_{{\cS}z(q)}(\bf D^+)$, hence $\cU$ comprises $r$ elements of
odd order, the same number of elements of order $2$ and $2r$
elements of order $4$. As before, Corollary \ref{HilDLS} allows us
to compute $g_{\cU}$ by means of the Hurwitz genus formula. The
case of $q-2q_0+1$ can be treated in a similar way. Therefore, we
have the following results.
\begin{proposition}
Let $\cU$ be a subgroup of order $4r$, with a divisor $r>1$ of
$q\pm 2q_0+1$. Then the quotient curve $\cX_\cU$ of $\cX$
associated to $\cU$ has genus
$$g_{\cU} \ = \left \{ \begin{array}{lll}
\frac{1}{4}[\frac{q_0(q-1)-1}{r}-(q_0-1)]  & \mbox{\em for \
$r\mid q+2q_0+1$\/}, \\ \frac{1}{4}[\frac{q_0(q-1)+1}{r}-(q_0+1)]
& \mbox{\em for \ $r\mid q-2q_0+1$\/}. \\
                         \end{array}
          \right.
$$
\end{proposition}

\section{Quotient curves arising from subgroups isomorphic to $\cS z({\bar q})$}

In this section we assume that $\bar{q}:=2^{2\bar{s}+1}$, with
$\bar{s}$ divisor of $s$ such that $2\bar{s}+1$ divides $2s+1$.
This is the arithmetical condition in order that ${\cS} z(q)$
contains a subgroup isomorphic to $\cS z({\bar q})$. Hence there
exists a subgroup $\cU$ of $\Aut(\bfq(\cX))$ isomorphic to $\cS
z({\bar q})$.
\begin{theorem}
\label{piccolosuzgen} Let $\bar{q}:=2^{2\bar{s}+1}$, with
$\bar{s}$ divisor of $s$ such that $2\bar{s}+1$ divides $2s+1$.The
quotient curve $\cX_{\cU}$ of $\cX$ associated to a subgroup $\cU$
isomorphic to ${\cS} z(\bar{q})$ has genus
$$g_{\cU}=\frac{q_0(q-1)-1+({\bar{q}}^2+1){\bar{q}}^2({\bar{q}}-1)+\Delta}
{({\bar{q}}^2+1){\bar{q}}^2({\bar{q}}-1)}$$
where
$$\Delta:=(\bar{q}^2+1)[(2q_0+2)(\bar{q}-1)+2\bar{q}(\bar{q}-1)]+
\bar{q}^2(\bar{q}^2+1)(\bar{q}-2)+\bar{q}^2(\bar{q}+2\bar{q}_0+1)
(\bar{q}-1)(\bar{q}-2\bar{q}_0).$$
\end{theorem}
\begin{proof} The group
$\cU$ has $(\bar{q}^2+1)(\bar{q}-1)$ elements of order $2$, and
$(\bar{q}^2+1)(\bar{q}^2-\bar{q})$ elements of order $4$.
Furthermore, $\cU$ has $\frac{1}{2}\bar{q}^2(\bar{q}^2+1)$
subgroups of order $\bar{q}-1$. Also, $\cU$ has
$\frac{1}{4}\bar{q}^2(\bar{q}+2\bar{q}_0+1)(\bar{q}-1)$ subgroups
of order $\bar{q}-2\bar{q}_0+1$, Finally, $\cU$ has
$\frac{1}{4}\bar{q}^2(\bar{q}-2\bar{q}_0+1)(q-1)$ subgroups of
order $\bar{q}+2\bar{q}_0+1$ By Corollary \ref{HilDLS}, $\deg
(\Diff (\cX|\cX_{\cU}))$ equals
$$(\bar{q}^2+1)[(2q_0+2)(\bar{q}-1)+2\bar{q}(\bar{q}-1)]+
\bar{q}^2(\bar{q}^2+1)(\bar{q}-2)+\bar{q}^2(\bar{q}+2\bar{q}_0+1)
(\bar{q}-1)(\bar{q}-2\bar{q}_0),$$ whence the assertion follows by
the Hurwitz genus formula.
\end{proof}
     
    \end{document}